\documentclass[6pt]{article}

 \normalbaselines
\newtheorem{Theorem}{Theorem}[section]

\newtheorem{Proposition}[Theorem]{Proposition}
\newtheorem{Lemma}[Theorem]{Lemma}

\newtheorem{Remark}[Theorem]{Remark}

\newcommand{\Proof}{{\em Proof }}

\usepackage[final]{graphicx}
\usepackage{amssymb}
\usepackage{amsmath}
\usepackage{pbox}
\usepackage{xcolor}
\usepackage{tcolorbox}
\usepackage{graphics}

\begin{document}

\title{ Volatility of volatility estimation: \\
central limit theorems for the Fourier transform estimator \\
and empirical study of the daily time series stylized facts
}

 \author{Giacomo Toscano\footnote{University of Firenze, Dept. of Economics and management, Firenze, Italy (corresponding author)}, Giulia Livieri\footnote{Scuola Normale Superiore, Pisa, Italy}, Maria Elvira Mancino\footnote{University of Firenze, Dept. of Economics and management, Firenze, Italy}, Stefano Marmi\footnote{Scuola Normale Superiore, Pisa, Italy}}

\date{May 25, 2022}

\maketitle

\textwidth=160 mm \textheight=250mm \parindent=8mm \frenchspacing
\vspace{ 3 mm}

\newpage

\begin{abstract}
{  We study the asymptotic normality of two {\color{black} feasible} estimators of the integrated volatility of volatility based on the Fourier methodology, which does not require the pre-estimation of the spot volatility.
We show that the bias-corrected estimator reaches the optimal rate \textcolor{black}{$n^{1/4}$}, while the estimator without bias-correction
 has a {slower} convergence rate and a smaller asymptotic variance.  Additionally, we provide simulation results {that} support the theoretical asymptotic distribution of the rate-efficient estimator and show the accuracy of the latter in comparison with a rate-optimal estimator based on the pre-estimation of the spot volatility.  Finally, using the rate-optimal Fourier estimator, we reconstruct the time series of the daily volatility of volatility of the S\&P500 and EUROSTOXX50 indices over long samples   and provide novel insight into the existence of stylized facts about the volatility of volatility  dynamics.}
\end{abstract}

\noindent {\bf Keywords:}   volatility of volatility, non-parametric estimation, central limit theorem, stochastic volatility,
Fourier analysis.

\noindent  {\bf JEL Classification:} C14, C58.

\section{Introduction}

In the last decades, different stochastic volatility models have been proposed to describe the evolution of asset prices, motivated by empirical studies on the patterns of volatilities in financial time series. Further, the availability of high-frequency data
has given impulse to devise statistical techniques aimed at the efficient estimation of model parameters in the stochastic volatility framework, e.g., the leverage and the volatility of volatility processes. {The estimation of  these model parameters is rather complicated,    the main difficulties being due to the fact that some factors are unobservable. In particular, the estimation of the volatility of volatility is a challenging task, because a pre-estimation of the spot volatility is typically required as a first step, due to the latency of the volatility process.}

Unlike the case of the integrated volatility, the non-parametric estimation of the integrated volatility of volatility is a relatively recent topic.
\cite{BNVe} propose a new class of stochastic volatility of volatility models, with an extra source of randomness, and show that the volatility of volatility
can be estimated non-parametrically by means of the quadratic variation of the preliminarily estimated squared volatility process, which they name {\sl pre-estimated spot variance based realized variance}. \cite{Vetter} proposes an estimator of the
integrated volatility of volatility  which is also based on increments of the pre-estimated spot volatility process and attains
the optimal convergence rate in the absence of noise. The common feature of these estimators is that they first reconstruct the unobservable volatility path via  some
consistent estimator thereof and then compute the volatility of volatility using the estimated paths as a proxy of the corresponding unknown paths.
The issue of estimating the volatility of volatility in the presence of jumps is studied in \cite{CuTe}:
first, the authors combine jump robust estimators of the integrated variance and the Fourier-Fej\'{e}r inversion formula to get an estimator of the instantaneous volatility path; secondly, they use again jump robust estimators of the integrated
volatility, in which they plug the estimated path of the volatility process, to obtain an estimator of the volatility of volatility. In the same spirit of \cite{BNVe,Vetter},  \cite{LiLiuZhang} also propose an estimator of the integrated volatility of volatility by means of a pre-estimation of the spot volatility, but, in order to extend the study to the case when the observed price process contains jumps and microstructure noise, the authors adopt a threshold pre-averaging estimator of the volatility, following \cite{Jing}.

In this paper, we focus on the estimation of the integrated volatility of volatility via the Fourier estimation method by \cite{MM}, which does not require the pre-estimation of the spot volatility. An early application of the Fourier methodology to identify the parameters (volatility of volatility and leverage) of stochastic volatility models has been proposed by \cite{BaMa2010}, where the authors prove a consistency result for the estimator of both the integrated leverage and volatility of volatility in the absence of noise. In the presence of microstructure noise, \cite{SanfeliciCuratoMancino} study the finite-sample properties of the Fourier estimator of the volatility of volatility introduced in \cite{BaMa2010} and show its asymptotic  unbiasedness.
However, the convergence rate of the estimator is not established, not even in the absence of microstructure noise contamination.

In the present paper we fill this gap. Specifically, after proving that the Fourier estimator of the volatility of volatility by \cite{SanfeliciCuratoMancino} has a sub-optimal rate of convergence, we define its bias-corrected version and prove  that it reaches the optimal convergence rate \textcolor{black}{$n^{1/4}$}.  { We also show that the non-corrected estimator with  slower rate of convergence displays a smaller asymptotic error variance.} {\color{black} Further,  we provide feasible versions of the two CLT's that exploit the product formula for the Fourier coefficients of the volatility of volatility and the fourth power of  {the} volatility. The same property of the Fourier coefficients is used in \cite{LMM} for the estimation of the quarticity.}

These asymptotic results   are supported by a simulation exercise, where we also compare the finite-sample performance of the rate-efficient Fourier estimator with that {of the rate-efficient realized estimator based on the pre-estimation of the spot volatility by \cite{ASJ}.} {\color{black} The comparative study suggests that the Fourier estimator works quite well on the daily horizon, while the performance of the realized estimator appears to be not satisfactory.} This feature may be related to the fact that, differently from the other volatility of volatility  estimators, which rely on the pre-estimation of the instantaneous volatility path  via a numerical differentiation, the Fourier approach relies only on the reconstruction of integrated quantities, i.e., the Fourier coefficients of the volatility. As it was early observed in \cite{MM}, this is a peculiarity of the Fourier estimator that renders the proposed method easily implementable and computationally stable.

Finally, we present an empirical exercise where the Fourier estimator is applied to   obtain the daily time series of the volatility of volatility of the S\&P500 and EUROSTOXX50 indices over, resp., the periods May 1, 2007 - August 6, 2021 and June 29, 2005 - May 28, 2021. As a result, we obtain some novel insight into the empirical regularities that characterize the daily dynamics of the volatility of volatility, which - to the best of our knowledge - up to now had been scarcely explored in the literature. Specifically, we find that the daily volatility of volatility of both the indices spikes in correspondence of periods of financial turmoil (e.g., during the financial crisis of 2008 and the outbreak of the COVID pandemic in 2020). Additionally,  we also find that it is usually positively (resp., negatively) correlated with the volatility (resp., the asset return), but appears to be less persistent than the volatility.  Finally, we observe that its empirical distribution is satisfactorily approximated by a log-normal distribution in years characterized by higher financial stability, as it is the case for the volatility.

\textcolor{black}{This novel insight appears to be valuable in view of the relevance of the volatility of volatility for scholars and practitioners.   Indeed, on the one hand, market operators regularly ``trade" the volatility of many financial asset classes via quoted and O.T.C. volatility derivatives (e.g., variance swaps, VIX futures and VIX options), hence the importance of the availability of accurate estimates of the volatility of the ``traded" volatility.} {\color{black} On the other hand, the need for efficient estimates of the volatility of volatility arises also in a number of technical tasks, e.g., the calibration of stochastic volatility of volatility models (\cite{BNVe}, \cite{SanfeliciCuratoMancino}),  the estimation of the leverage coefficient (\cite{KX}, \cite{AFLWY}), the inference of future returns \cite{Bolls}) and spot volatilities
(\cite{MZ2009}). Furthermore,  \cite{BFR} have recently provided empirical support to the dependence between the volatility of volatility of equity assets and structural sources of risk related to firms' characteristics.}

The paper is organized as follows. Section \ref{CQ} contains the assumptions and definitions. Section \ref{CLT} states the central limit theorems, which are supported by the simulation study in Section \ref{simulation}. Finally, Section \ref{empirical} contains the empirical results and Section \ref{concl} concludes. The proofs are given in Appendix A, while Appendix B contains some auxiliary lemmas on the Fej\'{e}r and Dirichlet kernels.

\section{Volatility of volatility estimators: definition and assumptions}
\label{CQ}

This section presents the general non-parametric stochastic volatility model which will be considered throughout the paper and defines two estimators of the integrated volatility of volatility based on the Fourier estimation method introduced in \cite{MM,MM09}. The class considered includes most of the continuous stochastic volatility models commonly used in high-frequency finance and is assumed (to cite one among many others) in Chapter 8.3 of \cite{ASJ}.

We make the following assumptions.

\medskip
{\rm (A.I)} The log-price process $p$ and the variance process $v$ are continuous It\^o semimartingales {on $[0,T]$} satisfying the stochastic differential equations

\begin{equation*}
\begin{cases}
dp(t)= \sigma(t)\,dW_t + b(t)\,dt\\
dv(t)= \gamma(t) dZ_t+\beta(t) dt
\end{cases}
\end{equation*}

\noindent where $v:=\sigma^2$, while $W$ and $Z$ are Brownian motions on a filtered probability space
$(\Omega, ({\cal F}_t)_{t\in [0,T]}, P)$ satisfying the usual conditions, possibly  correlated  \textcolor{black}{(in this regard, note that it is not restrictive to assume a constant correlation $\rho$).}

\par
{\rm (A.II)} The processes $\sigma$, $b$, $\gamma$ and $\beta$ are
continuous adapted stochastic processes defined on the same probability space
$(\Omega, ({\cal F}_t)_{t\in [0,T]}, P)$, such that for any $p\geq 1$,
$$
E\left[\int_0^T \sigma^p(t)dt\right]< \infty \ \ , \ E\left[\int_0^T b^p(t) dt\right]<\infty \ \ , \ E\left[\int_0^T \gamma^p(t)dt\right]< \infty \ \ , \ E\left[\int_0^T a^p(t) dt\right]<\infty.
$$
The processes are specified in such a way that \textcolor{black}{the spot volatility and volatility of volatility, resp. $\sigma$ and $\gamma$,} are a.s. positive. 

\par
{\rm (A.III)}
The process  $\gamma$ is a continuous It\^o semimartingale, whose drift and diffusion processes are
continuous adapted stochastic processes defined on the same probability space
$(\Omega, ({\cal F}_t)_{t\in [0,T]}, P)$.
 
\medskip

The assumptions (A.I)-(A.II)-(A.III) are standard in the non-parametric setting and are considered, e.g., in \cite{BNVe,ASJ,CuTe,Vetter,LiLiuZhang}.

\bigskip

By changing of the origin of time and scaling the unit of time, one can always modify the time window $[0,T]$ to $[0,2\pi]$.
Suppose that the asset log-price $p$ is observed at discrete, irregularly-spaced points in time on the grid
$\{ 0=t_{0,n} \leq  \ldots t_{i,n} \ldots \leq t_{n,n}=2\pi \}$.
For simplicity, we  omit the second index $n$. Denote $\rho(n):= ~ \max_{0\leq h \leq n-1}|t_{h+1}-t_{h}|$ and suppose that $\rho(n) \to 0$ as $n\to \infty$.

Consider the following interpolation formula
$$
p_{n}(t):= ~ \sum_{i=0}^{n-1} p(t_i) I_{[t_i,t_{i+1}[}(t).
$$
For any integer $k$, $|k|\leq 2N$, {\color{black} the discretized version of the Fourier coefficient $c_k(dp)$ is denoted by}
\begin{equation}
\label{FOUMULT1} c_k(dp_{n}):= {1\over {2\pi}} \sum_{i=0}^{n-1} e^{-{\rm i}kt_i} \, (p(t_{i+1})-p(t_i)),
\end{equation}
{\color{black}where the symbol ${\rm i}$ is the imaginary unit $\sqrt{-1}$}. Further, for any $|k|\leq N$, define the convolution formula
\begin{equation}
\label{CONVMULT}
c_k(v_{n,N}):={2\pi\over {2N+1}} \sum_{|s|\leq N} c_{s}(dp_n)c_{k-s}(dp_n).
\end{equation}
In \cite{MM09} it is proved that (\ref{CONVMULT}) is a consistent estimator of the $k$-th Fourier coefficient of the volatility process\footnote{Hereinafter, we will follow the relevant econometric literature by using the term volatility as a synonym of variance, thus referring to $\sigma^2(t)$ as the volatility process. Similarly for the volatility of volatility. } and in \cite{BaMa2010,SanfeliciCuratoMancino} it is shown that it is possible to derive an estimator of the integrated volatility of volatility by exploiting only the knowledge of the Fourier coefficients in (\ref{CONVMULT}), without the need of the preliminary estimation of the instantaneous volatility. This feature characterizes the Fourier method for estimating the volatility of volatility.  In fact, as far as we know, all other existing methods rely on the pre-estimation of the spot volatility, see \cite{ASJ,CuTe,Vetter,LiLiuZhang}. In general, these methods entail the pre-estimation of the spot volatility (in the absence or presence of noise contamination) as a first step; then, as a second step, a quadratic variation approach (e.g., the realized volatility formula) is applied to the pre-estimated spot volatility trajectory.

The estimator of the integrated volatility of volatility, defined in \cite{SanfeliciCuratoMancino}, is given by $2\pi$ times
\begin{equation}
\label{volvolestimator1}
c_0(\gamma^2_{n,N,M}):= {2\pi \over {M+1}} \sum_{|k| \leq M} \left( 1-{|k|\over {M+1}}\right) k^2 c_k(v_{n,N})c_{-k}(v_{n,N}),
\end{equation}
as $c_0(\gamma^2_{n,N,M})$ is the estimator of $c_0(\gamma^2)={1\over {2\pi}} \int_0^{2\pi}\gamma^2(t) dt$.
\textcolor{black}{ \cite{SanfeliciCuratoMancino} show that the estimator (\ref{volvolestimator1}) is consistent under the assumptions  {\rm (A.I)}-{\rm (A.II)}-{\rm (A.III)} and the conditions $N/n \to 0$ and $M^4/N\to 0$, in the absence of microstructure noise (Theorem 3.2). Moreover, they show that the estimator (\ref{volvolestimator1}) is asymptotically unbiased in the presence of microstructure noise (Theorem 4.2). However, they do not establish the rate of convergence and the asymptotic normality.
Further, note that the conditions on $n, \, N$ and $M$ assumed in Theorem 3.2 by \cite{SanfeliciCuratoMancino}  are only sufficient for the consistency of the estimator (\ref{volvolestimator1}); indeed, these conditions are not sharp, due to the fact that the focus of the paper is not  on the rate of convergence, but on the finite-sample properties of the estimator in the presence of microstructure noise. In the present paper,  we show that the convergence rate of the estimator (\ref{volvolestimator1}) is not optimal (see Theorem \ref{ASYMPT1}). Moreover, in Theorem \ref{ASYMPT2} we provide the optimal choices of the cutting frequencies $N$ and $M$ for the estimator (\ref{volvolestimator1}), in the absence of microstructure noise. In particular, note that Theorem \ref{ASYMPT2} assumes that the convolution parameter $N$ is equal to the Nyquist frequency $n/2$.  In the presence of microstructure contaminations, instead,  one needs to choose $N$ much smaller than $n/2$ to filter out the noise present in price observations, see Remark 4.3 in \cite{SanfeliciCuratoMancino}.}

\medskip

In order to obtain an estimator with the optimal rate of convergence in the absence of microstructure noise, a bias correction is needed and thus we consider the estimator
\begin{equation}
\label{volvolestimator3}
\widehat\gamma^2_{n,N,M}:= {2\pi \over {M+1}} \sum_{|k| \leq M}\left( 1-{|k|\over {M+1}}\right) k^2 c_k(v_{n,N})c_{-k}(v_{n,N}) - K \widehat \sigma^4_{n,N,M}
\end{equation}
where the constant $K$ is determined in (\ref{constantBiasCorrect}) and $\widehat \sigma^4_{n,N,M}$ is the Fourier estimator of the  quarticity, defined as
\begin{equation}
\label{ESTQUART} \widehat\sigma^4_{n,N,M}:= 2\pi \ \sum_{|k|\leq M} c_k(v_{n,N})c_{-k}(v_{n,N}).
\end{equation}
The asymptotic normality of the estimator (\ref{ESTQUART}) is studied in \cite{LMM}, while its properties in the presence of microstructure noise are studied in \cite{MSquart}.

Note that the estimator (\ref{volvolestimator3}) {differs from}  (\ref{volvolestimator1}) {for the presence of} the bias correction  $K \, \widehat \sigma^4_{n,N,M}$. This bias correction, while ensuring a faster rate of convergence, destroys the positivity of the estimator, see also \cite{BNHLS}. The estimator (\ref{volvolestimator1}) is instead positive.

\section{Central Limit Theorems}
\label{CLT}

In this section we study the asymptotic normality of the Fourier estimators of the integrated volatility of volatility defined by (\ref{volvolestimator1}) and
(\ref{volvolestimator3}) and prove that the  estimator (\ref{volvolestimator3}) reaches
the optimal rate of convergence \textcolor{black}{$n^{1/4}$},   at the cost of a de-biasing term, while the  estimator (\ref{volvolestimator1}) has a smaller asymptotic variance, at the cost of a slower convergence rate.

\begin{Theorem}
\label{ASYMPT1}
\par\noindent
Suppose that assumptions (A.I)-(A.II)-(A.III) hold. Let $N \rho(n) \sim c_N$ and $M\rho(n)^{1/2} \sim c_M$, where both $c_M$ and $c_N$ are positive constants\footnote{We stress the point that $c_N$ and $c_M$ are constants, i.e., they do not depend  on $M, N$.}. Then, as $n,N,M \to \infty$, the following stable convergence  in law holds:
$$\rho(n)^{-1/4}\left(\widehat\gamma^2_{n,N,M}- {1\over {2\pi}}\int_0^{2\pi}\gamma^2(t)dt\right)$$
$$\downarrow$$
$${\cal N}\left(0, {1\over {2\pi}}\int_0^{2\pi} \, K(c_M) \gamma^4(t) + K(c_M,c_N)\sigma^8(t)  + {\widetilde K}(c_M,c_N)\, \gamma^2(t) \, \sigma^4(t) \, dt \right),$$
where $K(c_M):={4\over 3} {1\over {c_M}}$, $K(c_M,c_N):={16\over {105}} \, c_M^3(1+2\eta(c_N/\pi))^2 $, ${\widetilde K}(c_M,c_N):={16\over{15}} \, c_M (1+2\eta(c_N/\pi))$ and
$\eta(a):=\frac{1}{2 a^2}r(a)(1-r(a))$, with
$r(a)=a-[a]$, being $[a]$ the integer part of $a$.
\end{Theorem}

Note that if $c_N=\pi$ or, equivalently, $N=n/2$ (i.e., the cutting frequency $N$ used for the estimation of the volatility coefficient given the log-prices is equal to the Nyquist frequency), then $\eta(c_N/\pi)=0$ and the asymptotic variance in Theorem \ref{ASYMPT1} becomes
\begin{equation}
\label{Var1}
{1\over {2\pi}}\int_0^{2\pi} {4\over 3}{1\over {c_M}}\gamma^4(t) + {16\over {105}} \, c_M^3 \sigma^8(t) + {16\over{15}} \, c_M
\gamma^2(t) \, \sigma^4(t) \, dt.
\end{equation}

\begin{Remark}
\label{comparison}
The realized volatility of volatility estimator (Th. 8.11 \cite{ASJ}) is obtained as the quadratic variation of the estimated spot volatility, with a de-biasing term depending on the quarticity. The underlying model is a continuous semimartingale  for the price, the volatility and the volatility of volatility.
The convergence rate of the estimator is {\color{black} $n^{1/4}$} and the asymptotic variance is
\begin{equation}
\label{Var2}
\int_0^T \, {151\over {70}} \beta \gamma^4(t)+{48\over {\beta^3}} \sigma^8(t) + {12\over \beta} \sigma^4(t) \gamma^2(t) \, dt.
\end{equation}
Letting $\beta= 1/c_M$,  the correspondence between the asymptotic variances (\ref{Var1}) and (\ref{Var2}) is easily seen, with the second and third terms   smaller in the case of the Fourier estimator. \textcolor{black}{Note that the estimator in \cite{ASJ} corresponds to (\ref{volvolestimator3}) multiplied by $2\pi$.}
\par\noindent
A similar approach as \cite{ASJ} is considered in \cite{LiLiuZhang}, but it is extended to obtain a consistent estimator in the presence of noisy data. To this aim, the authors first build an estimator of the spot volatility by means of a pre-averaging method to get rid of the noise contamination, then they compute the realized variance from the spot volatility estimates to obtain an estimator of the integrated volatility of volatility. Finally, they also need to correct for the bias of the obtained estimator. The rate of convergence is $n^{1/8}$ in the presence of noise and $n^{1/4}$ without noise.
\end{Remark}

{\color{black}In order to obtain a feasible CLT from Theorem \ref{ASYMPT1}, a consistent estimator of the conditional variance is needed. We exploit again the Fourier methodology to build a consistent estimator of
$${1\over {2\pi}} \int_0^{2\pi} \Lambda(t)dt$$
with
\begin{equation}
\label{asymvar_feasible}
\Lambda(t):= K(c_M) \gamma^4(t) + K(c_M,c_N)\sigma^8(t)  + {\widetilde K}(c_M,c_N) \gamma^2(t) \, \sigma^4(t).
\end{equation}
The result is detailed in Proposition \ref{Feasi1}. The key ingredients are the following Remark \ref{coefficients} and
the product formula for the Fourier coefficients, as studied in \cite{LMM}.}

\begin{Remark}
\label{coefficients}
{\color{black} The estimation of the integrated volatility of volatility relies on the convolution product
$$
 {2\pi\over {M+1}} \sum_{|k| \leq M}\left( 1-{|k|\over {M+1}}\right) c_k(dv) c_{-k}(dv),
$$
which allows computing the $0$-th Fourier coefficient of the volatility of volatility process.
The result is trivially extended to consider any continuous bounded function $h$ as
\begin{equation}
\label{hcase}
 {2\pi\over {M+1}} \sum_{|k| \leq M}\left( 1-{|k|\over {M+1}}\right) c_k(dv) c_{-k}(h \,dv),
\end{equation}
which leads to an estimator of
$${1\over {2\pi}} \int_0^{2\pi} h(t) \gamma^2(t) dt.
$$
In particular, for $h(t):= e^{-{\rm i}kt}$, the convolution product (\ref{hcase}) provides a formula for estimating the $k$-th Fourier coefficient of the volatility of volatility process $\gamma^2(t)$ (see also \cite{ClGl} for the analogous result in the case of the multivariate Fourier volatility estimator).}
\end{Remark}

{\color{black} Based on Remark \ref{coefficients},  for any integer $k$, $|k|\leq 2M$, we define
\begin{equation}
\label{quarticityK}
c_k(\sigma^4_{n,N,M}):= \sum_{|h|\leq M} c_h(v_{n,N}) c_{k-h}(v_{n,N})
\end{equation}
and
\begin{equation}
\label{unbiasedK}
\bar c_k( \gamma^2_{n,N,M}):=  {2\pi\over {M+1}} \sum_{|h| \leq M}\left( 1-{|h|\over {M+1}}\right) \, h(h-k) \, c_h(v_{n,N}) c_{k-h}(v_{n,N}) -K \, 2\pi \ c_k(\sigma^4_{n,N,M}),
\end{equation}
where $K$ is computed in (\ref{constantBiasCorrect}).
They are, resp., consistent estimators of $c_k(\sigma^4)$ and $c_k(\gamma^2)$, for any integer $k$.
The following result holds.}
\begin{Proposition}
\label{Feasi1}
\par\noindent
{\color{black} Suppose that assumptions (A.I)-(A.II)-(A.III) hold. Let $N \rho(n) \sim c_N>0$, $M\rho(n)^{1/2} \sim c_M>0$ and $LM^{-1}\to 0$  as $n,N,M,L \to \infty$. Let $\Lambda(t)$ be defined as in (\ref{asymvar_feasible}) and, further, define
\begin{equation}
\label{AsympVar_Est}
\Lambda_{n,N,M,L}:=K(c_M) { V}^{(1)}_{n,N,M,L} + {\widetilde K}(c_M,c_N) { V}^{(2)}_{n,N,M,L}+  K(c_M, c_N)  {V}^{(3)}_{n,N,M,L},
\end{equation}
where
$${ V}^{(1)}_{n,N,M,L}:=\sum_{|k|\leq L} \bar c_k( \gamma^2_{n,N,M})\bar c_{-k}( \gamma^2_{n,N,M}),
$$
$$
{ V}^{(2)}_{n,N,M,L}:= \sum_{|k|\leq L} \bar c_k(\gamma^2_{n,N,M})  c_{-k}( \sigma^4_{n,N,M}),
$$
$$
{ V}^{(3)}_{n,N,M,L}=: \sum_{|k|\leq L} c_{k}( \sigma^4_{n,N,M})c_{-k}( \sigma^4_{n,N,M}),
$$
where, in turn, $\bar c_k(\gamma^2_{n,N,M})$ is defined in (\ref{unbiasedK}), $c_{k}( \sigma^4_{n,N,M})$ is defined in (\ref{quarticityK}), and the constants $K(c_M)$, ${\widetilde K}(c_M,c_N)$ and $K(c_M,c_N)$ are specified in Theorem \ref{ASYMPT1}. As $n,N,M,L \to \infty$, the following convergence in probability holds:
$$
\Lambda_{n,N,M,L} \to  {1\over {2\pi}} \int_0^{2\pi} \Lambda(t)dt.
$$}
\end{Proposition}
\textcolor{black}{Thus, we have the following feasible CLT.}
\begin{Theorem}
\label{ASYMPT1feasible}
\par\noindent
{\color{black}Suppose that assumptions (A.I)-(A.II)-(A.III) hold. Let $N \rho(n) \sim c_N>0$ and $M\rho(n)^{1/2} \sim c_M>0$ and $LM^{-1}\to 0$  as $n,N,M,L \to \infty$. Let $\Lambda_{n,N,M,L}$ be defined as in (\ref{AsympVar_Est}).
Then, as $n,N,M,L \to \infty$, the following stable convergence  in law holds:
$$\rho(n)^{-1/4}\, {\widehat\gamma^2_{n,N,M}- \displaystyle {1\over {2\pi}} \int_0^{2\pi}\gamma^2(t)dt\over {\sqrt{\Lambda_{n,N,M,L}}  }} \to
{\cal N}\left(0,1 \right).$$}
\end{Theorem}

It is possible to obtain an estimator of the volatility of volatility without a bias-correction term and a smaller asymptotic variance, but the rate of convergence is slower, precisely $n^{\iota/2}$, with $\iota/2 \in (0,1/5)$. The estimator is simply given by (\ref{volvolestimator1}) multiplied by $2\pi$ and the
following result holds.
\begin{Theorem}
\label{ASYMPT2}
\par\noindent
Suppose that Assumptions (A.I)-(A.II)-(A.III) hold. Let $N \rho(n) \sim c_N>0$ and $M\rho(n)^{\iota} \sim c_M>0$, where $\iota \in (0,2/5)$. Then, as $n,N,M \to \infty$, the following stable convergence  in law holds:
$$\rho(n)^{-\iota/2}\left( c_0(\gamma^2_{n,N,M}) - {1\over {2\pi}}\int_0^{2\pi}\gamma^2(t)dt\right)\rightarrow
{\cal N}\left(0, {1\over {2\pi}}\int_0^{2\pi} {4\over 3} {1\over {c_M}}\gamma^4(t) \, dt \right).$$
\end{Theorem}

In order to build a feasible \textcolor{black}{CLT}, it is enough to apply the same methodology as for Theorem \ref{Feasi1}. In particular, under the conditions $N \rho(n) \sim c_N$ and $M\rho(n)^{\iota} \sim c_M$, where $\iota \in (0,2/5)$, a consistent estimator of the asymptotic variance is given by
\begin{equation}
\label{feas2}
\Gamma_{n,N,M,L}:=
{4\over 3} {1\over {c_M}} \sum_{|k|\leq L} c_k(\gamma^2_{n,N,M})c_{-k}( \gamma^2_{n,N,M}),
\end{equation}
where
\begin{equation}
\label{biasedK}
c_k(\gamma^2_{n,N,M}):={2\pi\over {M+1}} \sum_{|h| \leq M}\left( 1-{|h|\over {M+1}}\right) \, h(h-k) \, c_h(v_{n,N}) c_{k-h}(v_{n,N}).
\end{equation}
Therefore, the following holds.
\begin{Theorem}
\label{ASYMPT2feasible}
\par\noindent
Suppose that Assumptions (A.I)-(A.II)-(A.III) hold. Let $N \rho(n) \sim c_N>0$ and $M\rho(n)^{\iota} \sim c_M>0$, where $\iota \in (0,2/5)$. Then, as $n,N,M,L \to \infty$, the following stable convergence  in law holds:
$$\rho(n)^{-\iota/2} { c_0(\gamma^2_{n,N,M}) - \displaystyle {1\over {2\pi}} \int_0^{2\pi}\gamma^2(t)dt \over {\sqrt{\Gamma_{n,N,M,L}}}}\rightarrow
{\cal N}(0, 1).$$
\end{Theorem}

\begin{Remark}
The result in Theorem \ref{ASYMPT2} is in line with \cite{CuTe}, Th. \textcolor{black}{3.14}, where an estimator without bias correction is considered. \textcolor{black}{However, the proposed estimator relies on a smooth function of the plug-in spot volatility, where the latter is estimated with the Fourier method. Therefore, it differs from our estimator, which is based on the convolution formula of the Fourier coefficients of the volatility process.}
\end{Remark}

\section{Simulation study}
\label{simulation}

In this section we present a simulation study of the finite-sample performance of the rate-efficient  estimator
(\ref{volvolestimator3}). The objective of the study is to provide support to the asymptotic result in Theorem \ref{ASYMPT1}, \textcolor{black}{offer insight into the optimal selection of the frequency $M$, assess the robustness of the performance of the estimator to irregular sampling schemes} and illustrate a comparison of its accuracy  with that of the rate-efficient realized estimator by \cite{ASJ}.

\subsection{Simulation design}

We simulated discrete observations from two parametric models which satisfy Assumptions (A.I)-(A.II)-(A.III). The first model that  we simulated is the Heston model (see \cite{heston}):
\begin{equation}\label{sv}
\begin{cases}
dp(t)=\left(\mu-\frac{1}{2}v(t) \right)dt +\sigma(t) dW_t \\
dv(t)=\theta\left(\alpha-v(t)\right)dt +  \sqrt{\gamma^2 v(t)}dZ_t,
\end{cases}
\end{equation}
where $v(t):=\sigma^2(t)$,  $\mu \in \mathbb{R}$, $\theta,\alpha, \gamma>0$, and $\rho$ denotes the correlation between the Brownian motions $W$ and $Z$. Under the Heston model, the volatility of volatility
is given by   $\gamma^2(t)=\gamma^2 v(t)$.

The second model that we simulated is the stochastic volatility of volatility  model that appears in \cite{BNVee} and \cite{SanfeliciCuratoMancino}. The model is as follows:
\begin{equation}\label{svv}
\begin{cases}
dp(t)=\left(\mu-\frac{1}{2}v(t) \right)dt + \sigma(t) dW_t \\
dv(t)=\theta\left(\alpha-v(t)\right)dt +  \gamma(t)dZ_t \\
d\gamma^2(t)  =\chi \left(\eta -\gamma^2(t)\right)dt + \xi   \gamma(t) dY_t,
\end{cases}
\end{equation}
where $Y$ is a Brownian motion  independent of $W$ and $Z$, and $\mu \in \mathbb{R}$, $\theta,\alpha,  \chi, \eta, \xi >0$.

The parameter vectors used for the simulations of the models (\ref{sv}) and  (\ref{svv}) were, resp.:
\begin{itemize}
\item[-] $(\mu,\theta,\alpha, \gamma,\rho, x(0),v(0))=(0.1,5,0.2,0.5,-0.8,1,0.2)$;
\item[-]$(\mu,\theta,\alpha, \chi,\eta, \xi, \rho, x(0),v(0), \gamma^2(0))=(0.1,5,0.2,7,0.1,0.8,-0.8,1,0.2,0.1)$.
\end{itemize}

 Note that the selection  of a negative $\rho$ reproduces the presence of leverage effects.   For each model, we simulated \textcolor{black}{$10^4$}  trajectories with horizon $T=1/252$, i.e., with horizon equal to one trading day, corresponding to 6.5 hours.
For each trajectory, observations were simulated  on the equally-spaced grid with mesh equal to $1$ second.

\subsection{Finite-sample performance}
\label{fsp}

We assessed the finite-sample performance of the estimator
(\ref{volvolestimator3}) for increasing values of the sample size $n$, in order to provide numerical support to Theorem \ref{ASYMPT1}.  Specifically,  we considered values of $\rho(n)=T/n$ ranging between $5$ minutes and $1$ second. For what concerns the frequency $N$, which is needed for the convolution formula (\ref{CONVMULT}), we set $N=\left[c_N \rho(n)^{-1}\right]$ and selected $c_N=T/2$. This selection yields the value of $N$ equal to the Nyquist frequency $[n/2]$ and allows obtaining the smallest variance of the asymptotic error, see (\ref{Var1}) in Section \ref{CLT}. As for the frequency $M$, we set $M=\left[ c_M \rho(n)^{-1/2} \right]$
and optimized the value of the constant $c_M$ based on the (unfeasible) numerical minimization of the mean squared error (MSE). In this regard, we found that the MSE-optimal value of  $c_M$ is  equal to, resp., $0.05$  and $0.07$  for the models (\ref{sv}) and (\ref{svv}) (see also Subsection \ref{feaspro}, where  a feasible   procedure to select $M$ is also discussed).

Table \ref{table:sv} illustrates the finite-sample performance of the estimator (\ref{volvolestimator3}) under the two different data-generating processes considered. Specifically, Table \ref{table:sv} illustrates the MSE and the bias  for the different values of the   sampling frequency   $\rho(n) $. As expected, the bias and MSE improve  as $n$ is increased  for both the data-generating processes considered, thereby providing numerical support to Theorem \ref{ASYMPT1}. In particular,
note that the performance of the estimator  is still  satisfactory  for $\rho(n)$ equal to 5 minutes, the sampling frequency typically used in the absence of noise with empirical data.    Additionally, it is worth mentioning that the estimator never produced negative volatility of volatility estimates in this simulation study.

\begin{table}[h!]
\centering
{\color{black} \begin{tabular}{c |  c c | c c}
\multicolumn{1}{  c|}{ } & \multicolumn{2}{ c| }{Heston model } & \multicolumn{2}{ c }{stochastic vol-of-vol model}
 \\
\hline \hline
 $\rho(n)$  & MSE & Bias & MSE & Bias\\
 \hline
 5 minutes	&	$	 1.474	\cdot	10^{-8}	$	&	$	-5.674	\cdot	10^{-6}	$	 	&	$	 3.886   	\cdot	10^{-8}	$	&	$	1.455  	\cdot	10^{-6}	$	 \\
1 minute	&	$	4.425	\cdot	10^{-9}	$	&	$	-4.142	\cdot	10^{-6}	$	&	$	 1.939	\cdot	10^{-8}	$	&	$	1.304	\cdot	10^{-6}	$	\\
30 seconds	&	$	2.913	\cdot	10^{-9}	$	&	$	-3.801	\cdot	10^{-6}	$	&	$	 1.327	\cdot	10^{-8}	$	&	$	1.269	\cdot	10^{-6}	$	\\
5 seconds	&	$	9.985	\cdot	10^{-10}	$	&	$	-2.932	\cdot	10^{-7}	$	&	 $	8.113	\cdot	10^{-9}	$	&	$	4.161	\cdot	10^{-7}	$	 \\
1 second	&	$	4.229	\cdot	10^{-10}	$	&	$	-1.833	\cdot	10^{-7}	$	&	 $	6.199	\cdot	10^{-9}	$	&	$	3.644	\cdot	10^{-7}	$	 \\								
\end{tabular}}
\caption{MSE and bias of the estimator (\ref{volvolestimator3}) in correspondence of different values of $\rho(n)$. The true values of the daily integrated volatility of volatility are equal, on average, to $1.985\cdot 10^{-4}$ and $3.957\cdot 10^{-4}$ for, resp., the models (\ref{sv}) and (\ref{svv}).}
\label{table:sv}
\end{table}

\begin{figure}[h!]
\centering
 \includegraphics[width=\textwidth]{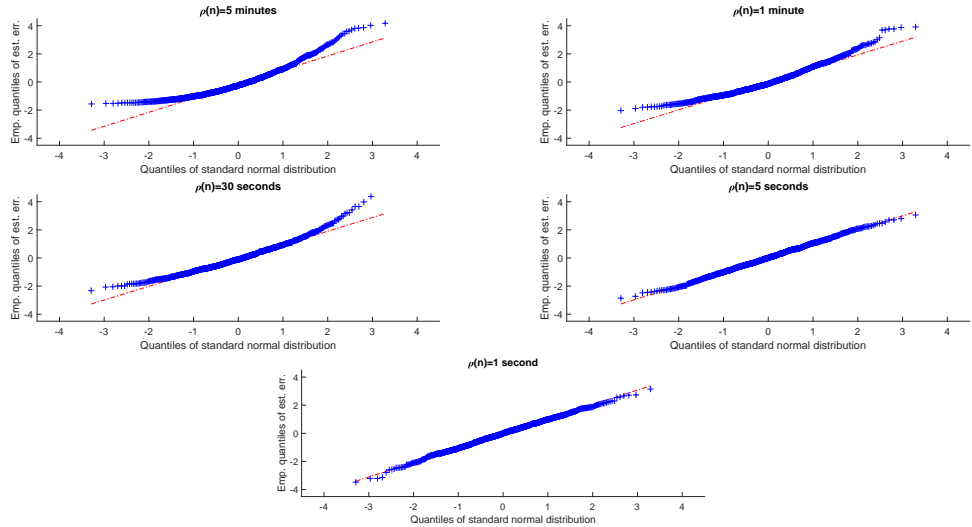}
 \caption{{\color{black}Comparison, for different values of $\rho(n)$, between the empirical quantiles of the standardized estimation error from Theorem \ref{ASYMPT1}, computed with simulated observations from the model (\ref{sv}), and the theoretical quantiles of a standard normal distribution.}}\label{cucuH}
\end{figure}

\begin{figure}[h!]
\centering
 \includegraphics[width=\textwidth]{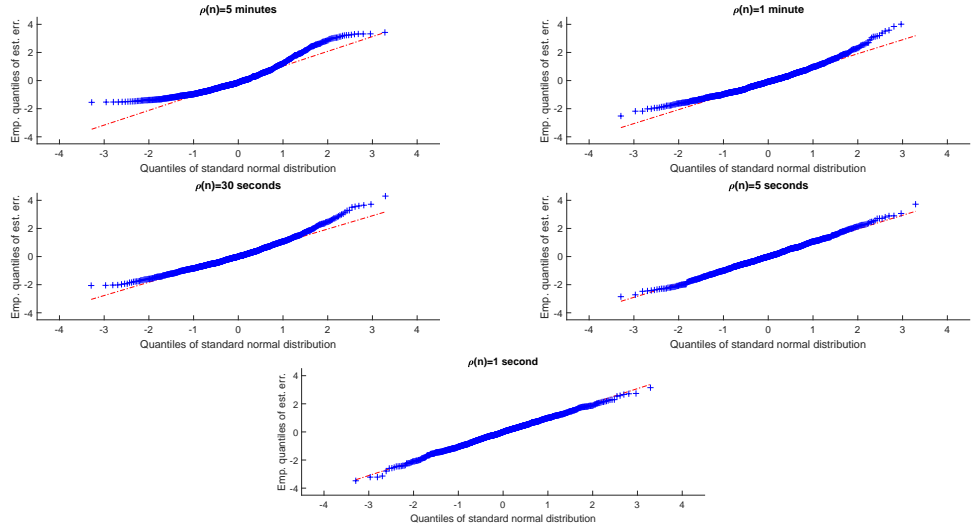}
 \caption{{\color{black}Comparison, for different values of $\rho(n)$, between the empirical quantiles of the standardized estimation error from Theorem \ref{ASYMPT1}, computed with simulated observations from the model (\ref{svv}), and the theoretical quantiles of a standard normal distribution.}}\label{cucuBNS}
\end{figure}

As additional support to the results in Thereom \ref{ASYMPT1}, the q-q plots in Figures  \ref{cucuH} and  \ref{cucuBNS}    offer a comparison between the empirical quantiles of the unfeasible standardized estimation error from Theorem \ref{ASYMPT1} and the theoretical quantiles of a standard normal distribution for different values of $\rho(n)$.  Both in the case of the model (\ref{sv})  and the model (\ref{svv}), as $\rho(n)$ becomes smaller, the approximation to the the standard normal distribution improves. In particular, while the approximation of the body of the distribution is satisfactory also for the largest $\rho(n)$ considered, i.e., $5$ minutes, the approximation in the tails becomes   accurate for $\rho(n)$ smaller or equal than $5$ seconds.

\subsection{Sensitivity to the frequency $M$}
\label{feaspro}

\begin{figure}[h!]
\centering
 \includegraphics[width=\textwidth]{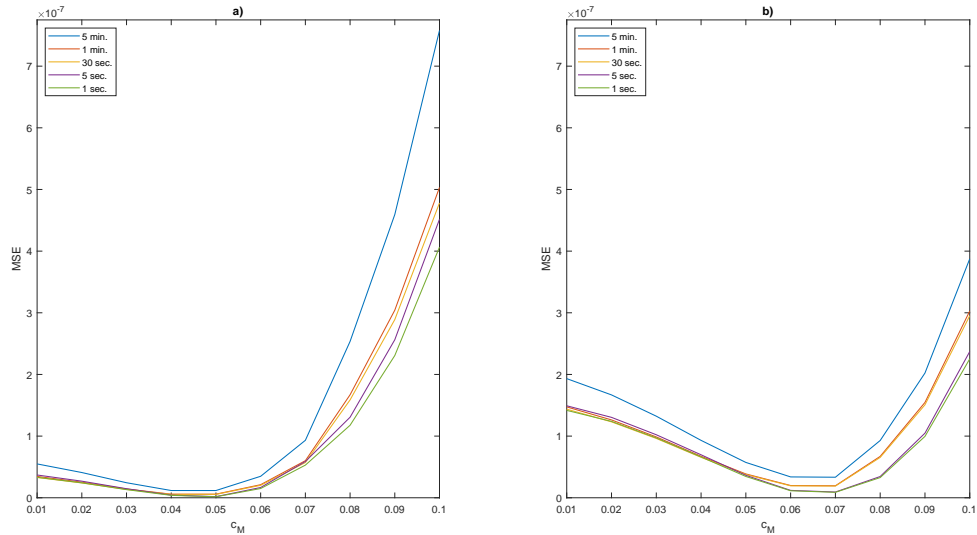}
 \caption{{\color{black} Sensitivity of the MSE of the estimator (\ref{volvolestimator3}) to different values of the constant $c_M$ in $M=[c_M\rho(n)^{-1/2}]$, for different values of $\rho(n)$. Panel a) refers to  the model (\ref{sv}), while panel b) refers to the model (\ref{svv}).}}\label{sensM}
\end{figure}

{\color{black} The careful selection of the frequency $M$, that is, the constant $c_M$, is key to efficiently implement the estimator (\ref{volvolestimator3}) with finite samples. Given the selection $M=\left[ c_M \rho(n)^{-1/2} \right]$, Figure \ref{sensM} shows the sensitivity of the MSE of the estimator to different values of the constant $c_M$ in the range $(0.01,1)$, for different values of $\rho(n)$. Based on Figure \ref{sensM}, the optimal MSE is achieved when $c_M$ is equal to, resp., $0.05$ and $0.07$ for the models (\ref{sv}) and (\ref{svv}), independently of $\rho(n)$. Further, it appears that the MSE is relatively flat for $c_M$ in the intervals $(0.04,0.06)$ and $(0.06,0.07)$. }

\textcolor{black}{As outlined in the following remark, it is possible to exploit the feasible Theorem \ref{ASYMPT1feasible} to optimize the selection of $c_M$ with empirical data.}

\begin{Remark}
{\color{black}The feasible selection of $c_M$  can be performed based on the adaptive procedure described in Section 2.4  of \cite{LiLiuZhang}. For the estimator (\ref{volvolestimator3}), the procedure can be summarized as follows. First, given an initial value for $c_M$, denoted by $c_{M,0}$, one computes the standard error of the estimator as $\rho(n)^{1/4}\sqrt{\Lambda_{n,N,M,L}}$, exploiting the feasible estimator of the asymptotic variance in (\ref{AsympVar_Est}). Then, the computation of the standard error is iterated in correspondence of increasing values of $c_M$ on a grid of step size $s$, that is, for $c_{M,j}, \, j\ge 1,$ such that $c_{M,j}-c_{M,j-1}=s \,\, \,  \forall j.$ Iterations stop when the absolute value of the marginal decrease in the standard error (w.r.t. the standard error recorded in correspondence to the initial value $c_{M,0}$) becomes smaller than a given threshold $\vartheta$. Unreported simulations show that, for $L=n^{1/4}$, the triple $(c_{M,0},s,\vartheta)=(0.03,0.01 ,0.25)$ provides an  estimation accuracy which is comparable  to that demonstrated in Subsection \ref{fsp} with the unfeasible optimization of the MSE. Moreover, the same simulations suggest that the procedure is fairly robust to   selections of the threshold   $\vartheta$ in the interval $(0.2,0.3)$.}
\end{Remark}

\subsection{Performance with irregular sampling}

{\color{black} So far, in the simulation study we assumed that prices were observable on the equally-spaced grid with mesh size $\rho(n)$ equal to $1$ second. However, the setup of Section \ref{CLT} allows for an irregular sampling scheme. To assess the robustness of the performance of the estimator (\ref{volvolestimator3}) to irregular sampling, we considered the case when observation times follow a  Poisson  process, that is,  durations between observations are drawn from an exponential distribution with mean $\lambda$ (see, e.g., \cite{book}, Chapter $3.3$).

Specifically, we considered three different values of $\lambda$, corresponding to an average duration $\delta$ of $1.25$, $1.5$ and $2$ seconds, and compared the resulting MSE and bias values with the case of regular sampling on the 1-second grid.} {\color{black} For the estimation with the Poisson scheme, we set $N=[n/2]$ and optimized  $M$  based on the minimization of the MSE. In this regard, we found that it is MSE-optimal to select a smaller $M$, compared to the regular-sampling case. Specifically, letting $M^*$ denote the optimal selection with regular 1-second sampling, numerical results suggest that it is optimal to select $M=[M^*/2]$. }

{\color{black} Table \ref{table:irr} reports the resulting bias and MSE, offering a comparison with the case where the a regular 1-second sampling scheme is adopted, corresponding to the last row of Table \ref{table:sv}. The results in Table \ref{table:irr} suggest that the Fourier estimator (\ref{volvolestimator3}) may still offer a satisfactory performance with irregular sampling schemes; in particular, it appears that the bias is relatively less affected than the MSE, compared to the regular-sampling case.}

\begin{table}[h!]
\centering
{\color{black}\begin{tabular}{c |  c c | c c}
\multicolumn{1}{  c|}{ } & \multicolumn{2}{ c| }{Heston model } & \multicolumn{2}{ c }{stochastic vol-of-vol model}
 \\
\hline \hline
{\color{black} $\delta$}  & MSE & Bias & MSE & Bias\\
 \hline
 {\color{black} $2$}	&	{\color{black}$	6.888	\cdot	10^{-10}	$}	&	{\color{black}$	-1.910	\cdot	10^{-7}	$}	&	{\color{black}$	9.672	\cdot	 10^{-9}	$}	&	{\color{black}$	3.840	\cdot	10^{-7}	$}	\\
{\color{black} $1.5$} 	&	{\color{black}$	 5.250	\cdot	10^{-10}	$}	&	{\color{black}$	-1.888	\cdot	10^{-7}	$}	&	{\color{black}$	7.564	\cdot	 10^{-9}	$}	&	{\color{black}$	3.776	\cdot	10^{-7}	$}	\\
{\color{black} $1.25$} 	&	{\color{black}$	4.662	\cdot	10^{-10}	$}	&	{\color{black}$	-1.865	\cdot	10^{-7}	$}	&	{\color{black}$	6.696	\cdot	 10^{-9}	$}	&	{\color{black}$	3.730	\cdot	10^{-7}	$}	\\ \hline
 {\color{black}regular 1-sec. sampling} 	&	$	4.229	\cdot	10^{-10}	$	&	$	-1.833	\cdot	10^{-7}	$	&	 $	6.199	\cdot	10^{-9}	$	&	$	 3.644	\cdot	10^{-7}	$	 \\								
\end{tabular}}
\caption{{\color{black} Comparison of MSE and bias of the estimator (\ref{volvolestimator3}) when observation times follow a Poisson scheme (with different values of the average duration $\delta$) with the case of regular 1-second sampling.}}  
\label{table:irr}
\end{table}

\subsection{Comparison with the performance of the rate-efficient realized estimator}
 \label{compa}

This subsection contains a comparative study of the finite-sample performance of the rate-efficient Fourier estimator (\ref{volvolestimator3})  and the rate-efficient realized estimator by \cite{ASJ} (see Remark \ref{comparison}).  We recall the definition of the latter. Let $\kappa(n)$ denote a sequence of integers such that $\kappa(n)\sim   \beta \rho(n)^{-1/2} $, $\beta>0, \, \rho(n):=T/n$. The estimator reads

\begin{equation}\label{aj}
\displaystyle \hat\gamma^2_{n }=\frac{3}{2 \kappa(n)}\sum_{i=1}^{n-2\kappa(n)+1}\Bigg(\Big( \widehat{ \sigma}_n^2( t_{i+k(n)}  ) - \widehat{\sigma}_n^2(t_{i} )\Big)^2 - \frac{4}{\kappa(n)}\big(\widehat{\sigma}_n^2(t_{i}\big))^2\Bigg),
\end{equation}
where
$$\displaystyle \widehat{\sigma}_{n }^2(t_i) = \frac{1}{\kappa(n)\rho(n)}\sum_{m=0}^{\kappa(n)-1} \Big(  p(t_{i+m}) - p(t_{i+m-1})\Big)^2$$
is the local estimator employed to pre-estimate the spot variance at time  $t_i=i{T/ n}$, $i=0,...,n$.   The  estimator  (\ref{aj}) is also studied in  \cite{Vetter}, where the author  replaces $\big(\widehat{\sigma}_n^2(t_{i}\big))^2$ with
 \begin{equation}\label{qvet}
  \displaystyle \widehat{\sigma}_n^4(t_i) = \frac{1}{3\kappa(n)\rho(n)^2 } \sum_{m=0}^{\kappa(n)-1} \Big(  p(t_{i+m}) - p(t_{i+m-1})\Big)^4.
 \end{equation}
The rate-efficient realized estimator considered in  \cite{Vetter} thus reads
\begin{equation}\label{vet}
 \hat{\bar\gamma}^2_{n}  =\frac{3}{2 \kappa(n)}\sum_{i=1}^{n-2\kappa(n)+1}\Bigg(\Big( \widehat{ \sigma}_n^2( t_{i+k(n)}  ) - \widehat{\sigma}_n^2(t_{i} )\Big)^2 - \frac{4}{\kappa(n)} \widehat{\sigma}_n^4(t_{i} ) \Bigg).
\end{equation}

{\color{black} For the comparison, we replicated the simulation study carried out in Subsection \ref{fsp}, this time using  the realized estimators (\ref{aj}) and (\ref{vet}) to obtain estimates of the daily integrated volatility of volatility. The implementation of realized estimators requires the selection of the tuning parameter $\beta$. After setting $\kappa(n)=[ \beta \rho(n)^{-1/2}]$, we selected  $\beta=0.04$ (resp., $\beta=0.06$) in the case of the model (\ref{sv}) (resp., \ref{svv}),  based on the unfeasible optimization of the MSE with 1-second samples.  Tables \ref{table:sv2} and \ref{table:sv3} summarize the results. By comparing  the latter with Table \ref{table:sv} in Subsection \ref{fsp}, it is immediate to see that   the performance of the realized estimators (\ref{aj}) and (\ref{vet}) is not satisfactory, both in terms of bias and MSE, compared to the case of the Fourier estimator (\ref{volvolestimator3}).  See also \cite{SanfeliciCuratoMancino} and \cite{TosRec} for similar considerations on the finite-sample performance of  realized volatility of volatility estimators. Moreover, note that the comparison for $\rho(n)$ equal to $5$ minutes is omitted, since the resulting bias and MSE of the realized estimators are larger than $1$ in absolute value.} Finally, simulations suggest that the use of the quarticity estimator in the de-biasing term in (\ref{vet}) does not improve the finite-sample performance.

\begin{table}[h!]
\centering
{\color{black} \begin{tabular}{c |  c c | c c}
\multicolumn{1}{  c|}{ } & \multicolumn{2}{ c| }{Heston model } & \multicolumn{2}{ c }{stochastic vol-of-vol model}
 \\
\hline \hline
$\rho(n)$  & MSE & Bias & MSE & Bias\\
 \hline
1 minute	&	{\color{black}$	1.800	\cdot	10^{-3}	$}	&	{\color{black}$	1.473	\cdot	10^{-2}	$}	&	{\color{black}$	1.655	\cdot	10^{-3}	$}	&	 {\color{black}$	1.299	\cdot	10^{-2}	$}	\\
30 seconds	&	{\color{black}$	 5.364	\cdot	10^{-4}	$}	&	{\color{black}$	1.388	\cdot	10^{-2}	$}	&	{\color{black}$	4.114	\cdot	10^{-4}	$}	&	 {\color{black}$	1.119	\cdot	10^{-2}	$}	\\
5 seconds	&	{\color{black}$	3.733	\cdot	10^{-4}	$}	&	{\color{black}$	 1.047	\cdot	10^{-2}	$}	&	 {\color{black}$	3.390	\cdot	10^{-4}	$}	&	 {\color{black}$	1.001	\cdot	10^{-2}	$}	 \\
1 second	&	{\color{black}$	3.322	\cdot	10^{-4}	$}	&	{\color{black}$	9.838	\cdot	10^{-3}	$}	&	{\color{black} $	2.999	\cdot	10^{-4}	$}	&	 {\color{black}$	9.555	\cdot	10^{-3}	$}	 \\								
\end{tabular}}
\caption{MSE and bias of the estimator (\ref{aj}) in correspondence of different values of $\rho(n)$. The true values of the daily integrated volatility of volatility are equal, on average, to $1.985\cdot 10^{-4}$ and $3.957\cdot 10^{-4}$ for, resp., the models (\ref{sv}) and (\ref{svv}).}
\label{table:sv2}
\end{table}

\begin{table}[h!]
\centering
{\color{black} \begin{tabular}{c |  c c | c c}
\multicolumn{1}{  c|}{ } & \multicolumn{2}{ c| }{Heston model } & \multicolumn{2}{ c }{stochastic vol-of-vol model}
 \\
\hline \hline
$\rho(n)$  & MSE & Bias & MSE & Bias\\
 \hline
1 minute	&	{\color{black}$	1.501	\cdot	10^{-3}	$}	&	{\color{black}$	1.575	\cdot	10^{-2}	$}	&	{\color{black}$	1.377	\cdot	10^{-3}	$}	&	 {\color{black}$	1.303	\cdot	10^{-2}	$}	\\
30 seconds	&	{\color{black}$	 5.461	\cdot	10^{-4}	$}	&	{\color{black}$	1.456	\cdot	10^{-2}	$}	&	{\color{black}$	4.336	\cdot	10^{-4}	$}	&	 {\color{black}$	1.122	\cdot	10^{-2}	$}	\\
5 seconds	&	{\color{black}$	3.783	\cdot	10^{-4}	$}	&	{\color{black}$	 1.091	\cdot	10^{-2}	$}	&	 {\color{black}$	3.302	\cdot	10^{-4}	$}	&	 {\color{black}$	9.998	\cdot	10^{-3}	$}	 \\
1 second	&	{\color{black}$	3.324	\cdot	10^{-4}	$}	&	{\color{black}$	9.840	\cdot	10^{-3}	$}	&	{\color{black} $	3.002	\cdot	10^{-4}	$}	&	 {\color{black}$	9.555	\cdot	10^{-3}	$}	 \\								
\end{tabular}}
\caption{MSE and bias of the estimator (\ref{vet}) in correspondence of different values of $\rho(n)$. The true values of the daily integrated volatility of volatility are equal, on average, to $1.985\cdot 10^{-4}$ and $3.957\cdot 10^{-4}$ for, resp., the models (\ref{sv}) and (\ref{svv}).}
\label{table:sv3}
\end{table}

\newpage

 \begin{Remark}
{\color{black} Unreported simulations show that the realized estimators (\ref{aj}) and (\ref{vet}) improve their finite-sample performance for larger values of the estimation horizon $T$. Specifically, realized estimators appear to achieve satisfactory accuracy, compared to the Fourier estimator (\ref{volvolestimator3}), when $T$ is equal to one year.  This is in line with the selection of $T$ equal to one year in the  numerical and  empirical high-frequency exercises  by \cite{LiLiuZhang}, where the performance of the noise- and jump-robust version of the realized estimator volatility of volatility estimator is investigated.} 
\end{Remark}

\section{Empirical analysis}
\label{empirical}

To the best of our knowledge, the empirical properties of the volatility of volatility of financial assets have been scarcely explored in the literature.
The aim of the empirical study presented in this section is thus to provide insight into the existence of stylized facts pertaining to the daily dynamics of the volatility of volatility.

In fact,  the numerical evidence presented in Section \ref{simulation} suggests that the Fourier methodology allows reconstructing the integrated volatility of volatility with satisfactory accuracy on daily intervals by means of the rate-efficient estimator (\ref{volvolestimator3}). Accordingly, in this section we use the estimator (\ref{volvolestimator3})
to obtain the daily volatility of volatility series for two market indices: the S\&P500 and the EUROSTOXX50. The periods considered for the analysis are, resp., May 1, 2007 - August 6, 2021 and June 29, 2005 - May 28, 2021.

\subsection{Data description, estimation and sample statistics}

For the empirical analysis we used the series of $5$-minute trade prices, recorded during trading hours. Specifically, for the S\&P500 index, we used the prices recorded between $9.30$ a.m. and $4$ p.m., while for the EUROSTOXX50 index we employed the prices recorded between $9$ a.m. and $5.30$ p.m. Days with early closure were discarded.

The estimation of the daily integrated volatility of volatility was performed  via the rate-efficient Fourier estimator (\ref{volvolestimator3}), without considering overnight returns. Before performing the estimation, we run the test by \cite{AX} on 5-minute series and found that the assumption of absence of noise could not be rejected at the $5\%$ significance level for both indices. Moreover, following \cite{WM14}, days with jumps were removed, based on the results of the test by \cite{LM08}, which was applied at the $1\%$ significance level. Overall, the number of days for which we estimated the volatility of volatility is $3343$ and $3522$ for, resp., the S\&P500 and EUROSTOXX50.

\begin{figure}[h!]
\centering
 \includegraphics[scale=0.35]{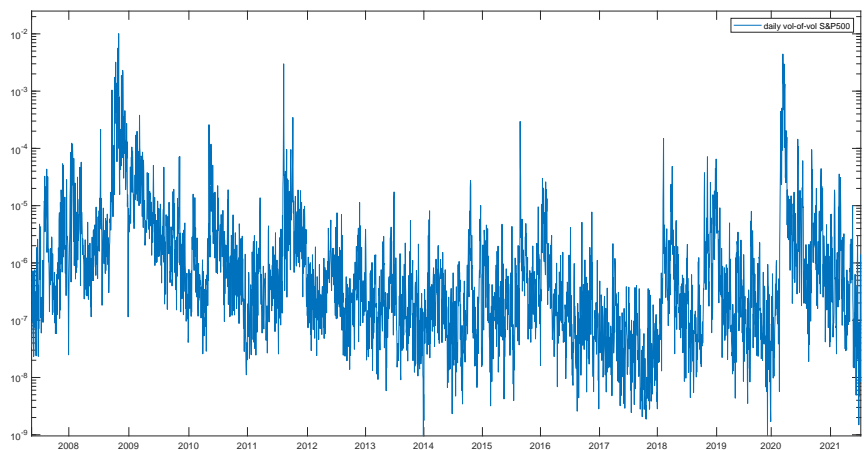}
 \caption{Daily volatility of volatility estimates for the S\&P500 index over the period  May 1, 2007- August 6, 2021. {A base-10 logarithmic scale is used for the y-axis.}}\label{usfig}
\end{figure}

\begin{figure}[h!]
\centering
 \includegraphics[scale=0.35]{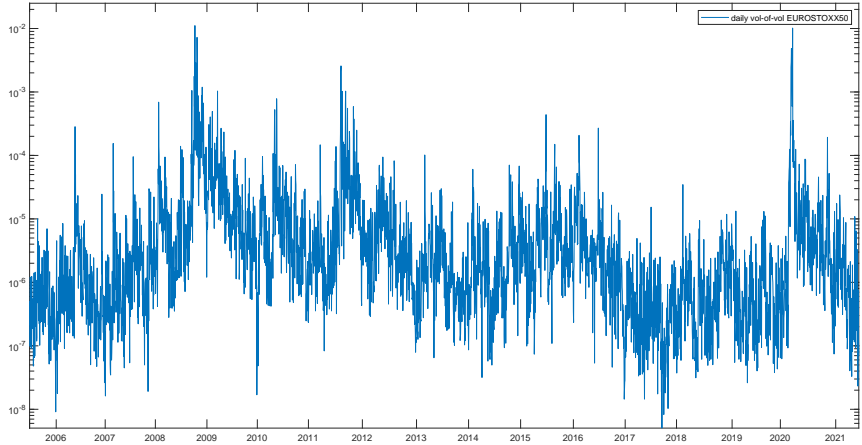}
 \caption{Daily volatility of volatility estimates for the EUROSTOXX50 index over the period June 29, 2005 - May 28, 2021. {A base-10 logarithmic scale is used for the y-axis.}}\label{eufig}
\end{figure}

\textcolor{black}{For the estimation, given the choice of the sampling frequency $\rho(n)=5$ minutes,   we set $N=[c_N\rho(n)^{-1}]$ and $M=[c_M\rho(n)^{-1/2}]$. Then we selected $c_N$ such that $N$ equals the Nyquist frequency (see Subsection \ref{fsp}) and $c_M$ based on the feasible procedure illustrated in Subsection \ref{feaspro}. The resulting values of $c_M$  are equal, on average, to $0.054$ and $0.065$ for, resp., the S\&P500 and the EUROSTOXX50.}  Figures \ref{usfig} and \ref{eufig}
display  the reconstructed trajectories of the daily volatility of volatility of the two indices, while Table \ref{sst} compares sample statistics of  volatility and volatility of volatility estimates. Daily volatility estimates were obtained via the Fourier estimator by \cite{MM}, applied to 5-minute returns\footnote{All the analyses appearing in this section and involving the computation of the volatility were also performed using volatility estimates obtained via the 5-minute realized variance  and the final outcome was pretty much the same.}. We note that all volatility of volatility estimates obtained are strictly positive.


\vspace{0.5cm}

 \begin{table}[h!]
\centering
\begin{tabular}{ c|  c  c c    c c   c c}
daily series    & mean &median & st. dev. & min & max & skew. & kurt.\\

 \hline\hline

SPX    vol.  & $ 8.598 \cdot 10^{-5}$ & $3.098 \cdot 10^{-5}$
 &  $  2.152  \cdot 10^{-4}$ & $ 1.701  \cdot 10^{-6}$
 &  $ 3.209 \cdot 10^{-3}$ &   $7.951  $ &   $83.902 $ \\

SPX  vol. of vol.   & $  \textcolor{black}{2.351} \cdot 10^{-5}$ & $ \textcolor{black}{3.867} \cdot 10^{-7}$ &   $ \textcolor{black}{2.608} \cdot 10^{-4}$ & $ \textcolor{black}{9.545 \cdot 10^{-10}}$ & $ \textcolor{black}{1.009} \cdot 10^{-2}  $ & $\textcolor{black}{24.442} $  & $\textcolor{black}{780.377}  $ \\   \hline

ESTX  vol. &  $1.093  \cdot 10^{-4}$ & $5.848 \cdot 10^{-5}$
 &  $1.944  \cdot 10^{-4}$ & $2.845 \cdot 10^{-6}$
 &  $3. 916 \cdot 10^{-3}$ &   $8.198$ &   $104.168$ \\

ESTX   vol. of vol. & $\textcolor{black}{3.233}  \cdot 10^{-5}$ & $\textcolor{black}{1.841} \cdot 10^{-6}$ &   $\textcolor{black}{3.267} \cdot 10^{-4}$ & $\textcolor{black}{5.054 \cdot 10^{-9}}$ & $\textcolor{black}{1.630} \cdot 10^{-2}  $ & $\textcolor{black}{24.726}$  & $\textcolor{black}{727.501}$ \\

			\end{tabular}
\caption{Sample statistics (mean, median, standard deviation, minimum, maximum, skewness, kurtosis) of the estimated daily integrated volatility and volatility of volatility of the S\&P500 (SPX)   and the EUROSTOXX50 (ESTX) indices over, resp., the periods May 1, 2007- August 6, 2021 and June 29, 2005 - May 28, 2021. }
\label{sst}
\end{table}

Figures \ref{usfig} and \ref{eufig} show  that the volatility of volatility spikes during financial crises, while remaining rather low and  stable during tranquil periods. Indeed,  both daily series reach their three highest peaks in correspondence of, resp., the global financial crisis starting at the end of 2008, the instabilities of the Euro area in the second part of 2011 and the outbreak of the COVID pandemic in the first months of 2020.

Moreover, based on Table \ref{sst}, we make the following remarks. First, the volatility of volatility is on average smaller than the volatility  in the case of both indices.  Secondly, the volatility of volatility appears to be more volatile than the volatility itself, as it displays  larger sample standard deviations and   maxima for both the estimated series. Finally, the volatility of volatility  appears to be much more skewed and leptokurtic than the volatility for both the indices.

 An analysis of the empirical regularities displayed by the reconstructed daily series of the volatility of volatility of the S\&P500 and the EUROSTOXX50 is illustrated in the next subsection.

 \subsection{Insight into volatility-of-volatility stylized facts}

The literature on the stylized facts related to the behavior of the volatility of financial assets is very rich (see, for instance,  \cite{ABDE}, \cite{PE} and \cite{corsi}, among many others). These include, e.g., clustering,  long memory, mean-reversion, log-normality  and leverage effects.  Nowadays, the volatility can be regarded, in some sense, as a traded asset itself. In fact, it is possible to ``trade" the volatility of many financial asset classes via quoted and O.T.C. volatility derivatives (e.g., variance swaps, VIX futures and VIX options).  Therefore, it may be of interest to evaluate which typical features of the daily volatility   actually apply to the daily volatility of the volatility itself.

\begin{figure}[h!]
\centering
 \includegraphics[scale=0.35]{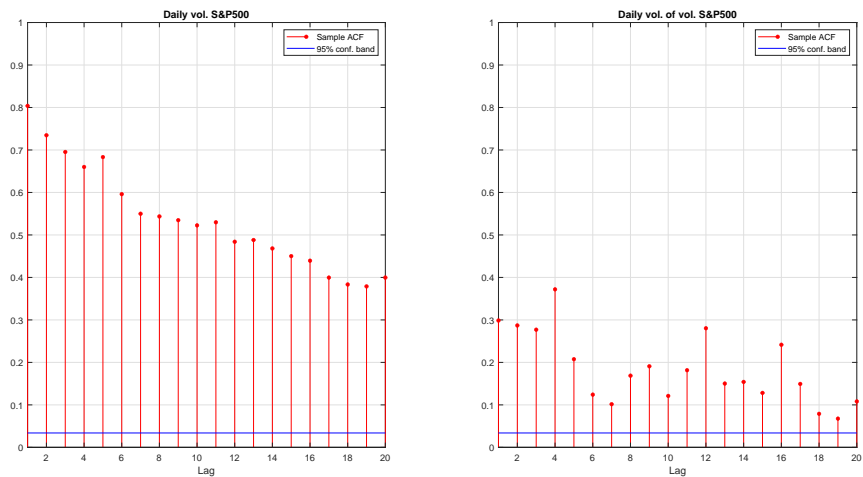}
 \caption{Sample autocorrelations (in red) of daily volatility (left panel) and volatility  of volatility (right panel)  for the S\&P500 index over the period  May 1, 2007- August 6, 2021. The $95\%$ confidence band (in blue) is computed under the null hypothesis of a Gaussian white noise process.}\label{acfus}
\end{figure}

\begin{figure}[h!]
\centering
 \includegraphics[scale=0.35]{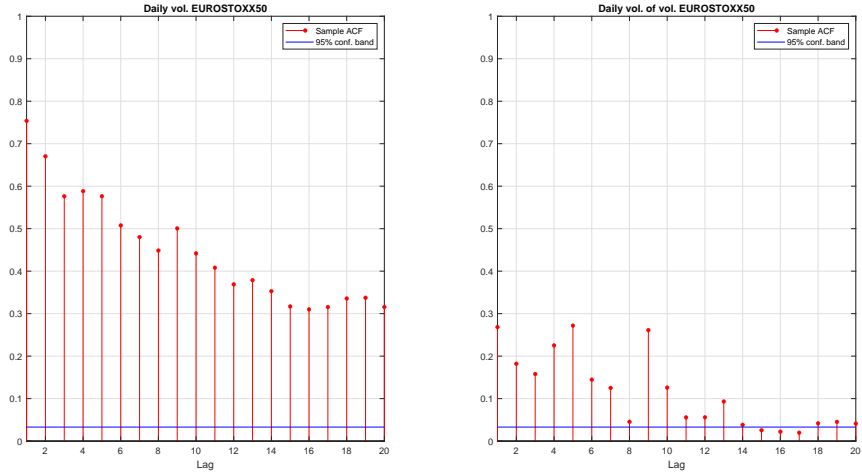}
 \caption{Sample autocorrelations (in red) of the daily volatility (left panel) and volatility  of volatility (right panel) for the EUROSTOXX50 index over the period June 29, 2005 - May 28, 2021. The $95\%$ confidence band (in blue) is computed under the null hypothesis of a Gaussian white noise process.}\label{acfeu}
\end{figure}

We have already observed in the previous subsection that the volatility of volatility shows clusters, being larger in correspondence of crises and smaller and less volatile during periods of economic stability. However, based on the observation of the sample auto-correlation function (see Figures \ref{acfus} and \ref{acfeu}), it appears to be less persistent than the volatility for both the indices considered.
As for the mean-reversion property, the Augmented Dickey-Fuller test rejects the hypothesis of a unit root for both volatility of volatility series analyzed, at the $0.01\%$ significance level.

We also examined the year-by-year correlation  of the daily volatility of volatility with, resp., the daily volatility and the daily log-return, computed as the difference between the closing and opening log-price. The dynamics of such correlations are summarized in Tables \ref{CORRSus} and \ref{CORRSeu}, where the values of  return-variance correlations, a rough proxy of the leverage effect, are also displayed for comparison\footnote{Note that the  correlations appearing in Tables \ref{CORRSus} and \ref{CORRSeu} are typically  pushed   towards zero by the presence of a finite-sample bias, see \cite{AFL}. For this reason, their true values are likely to be larger, in absolute value. However, obtaining unbiased and efficient estimates of these correlations goes beyond the scope of the exploratory  analysis proposed.}. For both the indices, we observe that the yearly correlation between the log-return and the volatility of volatility tends to be negative and to follow the return-variance correlation closely, although being most often smaller in absolute value. This result may suggest  the existence of a ``second-order" leverage effect: what we observe is in fact that when the asset price decreases, not only  the volatility increases, due the asset becoming riskier, but also the  volatility of volatility - which can be seen as a proxy of the  uncertainty about the amount of risk perceived by market operators, that is, the ``volatility of risk" - becomes larger.  The yearly correlation between the volatility of volatility and the volatility is instead    positive and close to {\color{black} $0.7$} on average for both the indices. This is consistent with the presence of volatility-of-volatility peaks in periods of higher market volatility.

 \begin{table}[h!]
\centering
\begin{tabular}{ c|  c  c c    }
 year	&	vol. of vol. - vol. 	&	vol.- ret.  	&	vol. of vol. - ret.  		\\
\hline	\hline							
2007	&	\textcolor{black}{	    0.730	}	&	-0.126	&	 \textcolor{black}{	-0.030}	\\	
2008	&	\textcolor{black}{	    0.689	}	&	0.049	&	\textcolor{black}{	    0.054	} 	\\	
2009	&	\textcolor{black}{	    0.727	}	&	-0.063	&	\textcolor{black}{	    0.067	} 	\\	
2010	&	\textcolor{black}{	    0.827	}	&	-0.199	&	\textcolor{black}{	   -0.073	} 	\\	
2011	&	\textcolor{black}{	    0.728	}	&	0.019	&	\textcolor{black}{	    0.213	} \\	
2012	&	\textcolor{black}{	    0.686	}	&	-0.231	&	\textcolor{black}{	    0.022	} 	\\	
2013	&	\textcolor{black}{	    0.706	}	&	-0.354	&	\textcolor{black}{	   -0.129	} \\	
2014	&	\textcolor{black}{	    0.842	}	&	-0.140	&	\textcolor{black}{	    0.030	} 	\\	
2015	&	\textcolor{black}{	    0.582	}	&	-0.153	&	\textcolor{black}{	   -0.281	} \\	
2016	&	\textcolor{black}{	    0.849	}	&	-0.122	&	\textcolor{black}{	   -0.151	} 	\\	
2017	&	\textcolor{black}{	    0.638	}	&	-0.348	&	\textcolor{black}{	   -0.160	} \\	
2018	&	\textcolor{black}{	    0.758	}	&	-0.218	&	\textcolor{black}{	   -0.125	}\\	
2019	&	\textcolor{black}{	    0.858	}	&	-0.125	&	\textcolor{black}{	   -0.026	}	\\	
2020	&	\textcolor{black}{	    0.821	}	&	0.090	&	\textcolor{black}{	    0.172	} 	\\	
2021	&	\textcolor{black}{	    0.913	}	&	-0.305	&	\textcolor{black}{	   -0.219 	} 	\\	
\hline								
average	&	\textcolor{black}{0.757}	&	-0.147	&	\textcolor{black}{-0.043}

			\end{tabular}
\caption{S\&P500 index: sample yearly correlations of the daily volatility  of volatility with the corresponding daily volatility and daily log-return. }
\label{CORRSus}
\end{table}

 \begin{table}[h!]
\centering
\begin{tabular}{ c|  c  c c    }
 year	&	vol. of vol. - vol. 	&	vol.- ret.  	&	vol. of vol. - ret.  		\\
\hline	\hline	

2005	&	\textcolor{black}{	0.6906	}	&	-0.1463	&	\textcolor{black}{		-0.1740	}	\\
2006	&	\textcolor{black}{	    0.593	}	&	   -0.119	&	\textcolor{black}{		   -0.101	}	\\
2007	&	\textcolor{black}{	    0.689	}	&	   -0.237	&	\textcolor{black}{		   -0.252	}	\\
2008	&	\textcolor{black}{	    0.833	}	&	   -0.130	&	\textcolor{black}{		   -0.053	}	\\
2009	&	\textcolor{black}{	    0.621	}	&	   -0.144	&	\textcolor{black}{		   -0.009	}	\\
2010	&	\textcolor{black}{	    0.804	}	&	   -0.107	&	\textcolor{black}{		   -0.040	}	\\
2011	&	\textcolor{black}{	    0.864	}	&	   -0.072	&	\textcolor{black}{		   -0.156	}	\\
2012	&	\textcolor{black}{	    0.718	}	&	   -0.035	&	\textcolor{black}{		    0.023	}	\\
2013	&	\textcolor{black}{	    0.685	}	&	   -0.141	&	\textcolor{black}{		   -0.052	}	\\
2014	&	\textcolor{black}{	    0.742	}	&	   -0.128	&	\textcolor{black}{		   -0.074	}	\\
2015	&	\textcolor{black}{	    0.584	}	&	   -0.111	&	\textcolor{black}{		    0.052	}	\\
2016	&	\textcolor{black}{	    0.732	}	&	   -0.327	&	\textcolor{black}{		   -0.275	}	\\
2017	&	\textcolor{black}{	    0.552	}	&	   -0.172	&	\textcolor{black}{		   -0.196	}	\\
2018	&	\textcolor{black}{	    0.791	}	&	   -0.327	&	\textcolor{black}{		   -0.303	}	\\
2019	&	\textcolor{black}{	    0.743	}	&	   -0.180	&	\textcolor{black}{		   -0.016	}	\\
2020	&	\textcolor{black}{	    0.705	}	&	   -0.276	&	\textcolor{black}{		   -0.084	}	\\
2021	&	\textcolor{black}{	    0.762	}	&	   -0.185	&	\textcolor{black}{		   -0.081	}	\\

\hline								
average	&	\textcolor{black}{	    0.712	}	&	-0.167	&	\textcolor{black}{	   -0.105	}		
			\end{tabular}
\caption{EUROSTOXX50 index: sample yearly correlations of the daily volatility  of volatility with the corresponding daily volatility and daily log-return.  }
\label{CORRSeu}
\end{table}

Finally, we test for the Gaussianity of the logarithmic volatility and volatility of volatility estimates, using the Jarque-Bera and Anderson-Darling tests at the $5\%$ significance level. The years in which both tests reject the null hypothesis of Gaussianity are 2008, 2011, 2016 and 2020, that is, the years in the sample that were the most characterized by market turmoil (in order: the global financial crisis, the Euro-area instability phase, Brexit and the outbreak of the COVID pandemic). This happens for both the quantities tested and both the indices analyzed, thus suggesting that the log-normal approximation for the distribution of the volatility  and the volatility of volatility is more satisfactory in periods of market stability.

\section{Conclusions}
\label{concl}
This paper fills a gap in the literature on financial econometrics by deriving the convergence rate of the Fourier estimator of the volatility of volatility.  In this regard, we showed that the bias-corrected version of the estimator reaches the optimal rate \textcolor{black}{$n^{1/4}$}, while the estimator without bias-correction achieves a sub-optimal rate, but has a smaller asymptotic variance.

Further, we presented a numerical study that shows   that the rate-optimal Fourier estimator of the volatility of volatility performs well in finite samples, even at the relatively small daily estimation horizon, where the competing rate-efficient realized estimator shows a poor performance.

Finally, we applied the Fourier estimator to multi-year samples of  S\&P500 and EUROSTOXX50 observations    and gained some new knowledge   about the empirical regularities that characterize the daily dynamics of the volatility of volatility, a topic which so far had been  scarcely explored in the literature.

\bigskip

\section{Appendix A: proofs}
\label{Proofs}

{\color{black}The proofs of Theorems \ref{ASYMPT1}, \ref{ASYMPT1feasible}, \ref{ASYMPT2} and \ref{ASYMPT2feasible} are illustrated in the next subsections}. Some preliminary remarks are useful.

\begin{Remark}
\label{regularity}
As every continuous process is locally bounded, all processes appearing in the proofs are locally bounded. Moreover, standard localization procedures (see, e.g., \cite{ASJ}) allow assuming that any locally bounded process is actually bounded and almost-surely positive processes can be considered as bounded away from zero.
\end{Remark}

\begin{Remark}
\label{boundary}
In \cite{MM09}, Lemma 2.2, it is proved that the drift component of the semi-martingale model gives no contribution to the convolution formula (\ref{CONVMULT}). Therefore, we will refer to the drift-less model in Assumption (A.I). Moreover,
as observed in \cite{MM}, we can assume that $p(0)=p(2\pi)$ and $v(0)=v(2\pi)$. In fact, if $p(0)\not= p(2\pi)$ (similarly for the process $v$), we introduce
$$\widetilde p(t)= p(t)- {t\over {2\pi}}[p(2\pi)-p(0)].$$
Then, while $\widetilde p$ satisfies the required assumption, at the same time the volatility and co-volatilities estimations are not affected by a modification of the drift as above.
From the point of view of the modeling, we may consider
$$d\widetilde p(t)= \sqrt{v(t)} dW(t),$$
$$d \widetilde v(t) = \gamma(t) \, dZ(t),$$
being $v(t)= \widetilde v(t) - {t\over {2\pi}}[{\widetilde v}(2\pi)-{\widetilde v}(0)].$ In fact, for any $k\not= 0$, it holds $c_k(dv)=c_k(d\widetilde v)$, while the $0$-th Fourier coefficient, $c_0(dv)$, is not contributing to the definition of the Fourier  estimator of the volatility of volatility.
The situation would change in the case when one wishes to estimate the spot volatility.
However, this is not an issue of the present study, as the estimation of the spot volatility  is not required.
\end{Remark}

\subsection{Preliminary decomposition}
\label{PrelDecom}

Given the discrete time observations
$\{ 0=t_{0} \leq  \ldots \leq t_{i} \leq \ldots \leq t_{n}=2\pi \}$, we use the notation in continuous time by letting $\varphi_n(t):= \sup \{t_j : t_j \leq t\}$, for the sake of simplicity.
From the It\^o formula we have the following decomposition of the term (\ref{CONVMULT})
\begin{equation}
\label{coeffVolDecomp}
c_{k}(v_{n,N}):= A_{k,n}+B_{k,n,N}+ C_{k,n,N},
\end{equation}
where
\begin{equation}
\label{coeffVolDecompA}
A_{k,n}:={1\over {2\pi}} \int_0^{2\pi} e^{-ik\varphi_n(s)} v(s) ds,
\end{equation}
\begin{equation}
\label{coeffVolDecompB}
 B_{k,n,N}:=
 {1\over {2\pi}} \int_0^{2\pi} e^{-ik\varphi_n(s)} \int_0^s D_{N}(\varphi_n(s)-\varphi_n(u)) \sigma(u) dW_u \, \sigma(s) dW_s,
\end{equation}
\begin{equation}
\label{coeffVolDecompC}
C_{k,n,N}:= {1\over {2\pi}} \int_0^{2\pi} \int_0^s e^{-ik\varphi_n(u)} D_{N}(\varphi_n(s)-\varphi_n(u)) \sigma(u) dW_u \, \sigma(s) dW_s.
\end{equation}
It follows that {(\ref{volvolestimator1})} is equal to the following sum:
\begin{equation}
\label{Decomposition0}
{2\pi\over {M+1}}\sum_{|k| \leq M} \left(1-{|k|\over {M+1}}\right) k^2 \, A_{k,n}A_{-k,n}
\end{equation}
\begin{equation}
\label{Decomposition1}
+{2\pi\over {M+1}} \sum_{|k| \leq M} \left(1-{|k|\over {M+1}}\right) \, k^2 \, 2(A_{k,n}B_{-k,n,N}+ A_{k,n}C_{-k,n,N})
\end{equation}
\begin{equation}
\label{Decomposition2}
+{2\pi\over {M+1}}\sum_{|k| \leq M} \left(1-{|k|\over {M+1}}\right) \, k^2 \,  (B_{k,n,N}B_{-k,n,N}+ 2B_{k,n,N}C_{-k,n,N}+ C_{k,n,N}C_{-k,n,N}).
\end{equation}

Throughout the proofs, we denote 
\begin{equation}
\label{integradirichlet}
Y_{n,N}(t,s):=\int^t_0 D_{N}(\varphi_n(s)-\varphi_n(u))\sigma(u)dW_u
\end{equation}
and, for brevity, we will use the following notation for the Dirichlet and the Fej\'{e}r kernels (see also the Appendix B):
\begin{equation}
\label{notazioneDF}
D_{N,n}(s-u):=D_{N}(\varphi_n(s)-\varphi_n(u)), \ \ \ \ F_{M,n}(s-u):= F_M(\varphi_n(s)-\varphi_n(u)).
\end{equation}
Similarly, for the derivatives of the Fej\'{e}r kernel\footnote{We stress the point that the notation $D_{N,n}(s-u)$ refers to $D_{N}(\varphi_n(s)-\varphi_n(u))$, which is a function of two time variables, $(s,u)$, as stated in (\ref{notazioneDF}). It is not to identify it with $D_{N}(\varphi_n(s-u))$. Notation (\ref{notazioneDF}) is maintained in order to highlight the role of the convolution product, which is the main characteristic of the Fourier estimation method in \cite{MM09}.} we will use $F^{\prime}_{M,n}$ and $F^{\prime\prime}_{M,n}$.

In order to identify the different contribution of all terms, we start with equation (\ref{Decomposition0}), which can be written as
\begin{equation}
\label{continuous}
{2\pi\over {M+1}}\sum_{|k| \leq M} \left(1-{|k|\over {M+1}}\right) \, k^2 c_{k}(v)c_{-k}(v)
\end{equation}
\begin{equation}
\label{AAnew}
+{2\pi\over {M+1}}\sum_{|k| \leq M} \left(1-{|k|\over {M+1}}\right) k^2 \, \left(A_{k,n}A_{-k,n}-  c_{k}(v)c_{-k}(v)\right).
\end{equation}
We consider now the term (\ref{Decomposition1}). An application of the It\^o formula shows that it is equal to
$$ 2(AB^{(i)}_{M,n,N}+AB^{(ii)}_{M,n,N}+AC^{(i)}_{M,n,N}+AC^{(ii)}_{M,n,N}),$$
where
$$
AB^{(i)}_{M,n,N}:= {1\over {2\pi}}\int_0^{2\pi} \int_0^s {1\over {M+1}} \sum_{|k|\leq M} \left(1-{|k|\over {M+1}}\right)k^2 e^{{\rm i}k(\varphi_n(s)-\varphi_n(u))}v(u) du \, Y_{n,N}(s,s) \sigma(s) dW_{s}
$$
\begin{equation}
\label{ABi2}
= - \int_0^{2\pi} {1\over {2\pi}} \int_0^s {1\over {M+1}} F_{M,n}^{\prime \prime}(s-u) v(u) du \, Y_{n,N}(s,s) \sigma(s) dW_s,
\end{equation}
$$
AB^{(ii)}_{M,n,N}:=
\int_0^{2\pi} {1\over {2\pi}}\int_0^s {1\over {M+1}}\sum_{|k|\leq M} \left(1-{|k|\over {M+1}}\right)k^2 e^{-{\rm i}k(\varphi_n(s)-\varphi_n(u))} \, Y_{n,N}(u,s) \sigma(u) dW_{u} \, v(s) ds
$$
\begin{equation}
\label{ABii2}
=-\int_0^{2\pi}{1\over {2\pi}} \int_0^s {1\over {M+1}}F_{M,n}^{\prime \prime}(s-u) \, Y_{n,N}(u,s) \, \sigma(u) dW_{u} \, v(s) ds,
\end{equation}
and, letting
$$Y_{k,n,N}(t,s):= \int_0^t e^{-{\rm i}k\varphi_n(u)} D_N(\varphi_n(s)-\varphi_n(u)) \sigma(u) dW_u,$$
then
$$
AC^{(i)}_{M,n,N}:={1\over {M+1}}  \sum_{|k|\leq M} \left(1- {|k|\over {M+1}}\right) k^2 {1\over {2\pi}}\int_0^{2\pi} \int_0^s e^{-{\rm i}k \varphi_n(u)}v(u) du \, Y_{-k,n,N}(s,s) \sigma(s) dW_s
$$
\begin{equation}
\label{ACia}
=- \int_0^{2\pi} {1\over {2\pi}}\int_0^s \int_0^u {1\over {M+1}} F_{M,n}^{\prime \prime}(u-u_1)D_{N,n}(s-u_1)\sigma(u_1)dW_{u_1} \, v(u) du \, \sigma(s) dW_s
\end{equation}
\begin{equation}
\label{ACib}
-  \int_0^{2\pi} {1\over {2\pi}}\int_0^s \int_0^u {1\over {M+1}} F_{M,n}^{\prime \prime}(u-u_1)v(u_1) du_1 D_{N,n}(s-u)\sigma(u)dW_{u} \, \sigma(s) dW_s,
\end{equation}
and
$$
AC^{(ii)}_{M,n,N}:=
{1\over {M+1}}  \sum_{|k|\leq M} \left(1- {|k|\over {M+1}}\right) k^2 {1\over {2\pi}}\int_0^{2\pi} \int_0^s  Y_{-k,n,N}(s_1,s) \sigma(s_1)dW_{s_1} \, e^{-{\rm i}k\varphi_n(s)} v(s) ds
$$
\begin{equation}
\label{ACii2}
=- \int_0^{2\pi}  {1\over {M+1}}{1\over {2\pi}}\int_0^s \int_0^{s^{\prime}} F_{M,n}^{\prime \prime}(s-u) D_{N,n}(s-u)\, \sigma(u)dW_u \, \sigma(s_1) dW_{s_1} \, v(s) ds.
\end{equation}

Finally, we consider the term (\ref{Decomposition2}). Using the fact that the Fej\'{e}r kernel is an even function, it can be re-written as
$$
{2\pi \over {M+1}}\sum_{|k| \leq M} \left(1-{|k|\over {M+1}}\right) \, k^2 \, (B_{k,n,N}B_{-k,n,N}+ 2B_{k,n,N}C_{-k,n,N}+ C_{k,n,N}C_{-k,n,N})
$$
$$
=BB_{M,n,N}^{(i)}+2 BB_{M,n,N}^{(ii)}+ 2 (BC_{M,n,N}^{(i)}+BC_{M,n,N}^{(ii)}+BC_{M,n,N}^{(iii)})+CC_{M,n,N}^{(i)}+2CC_{M,n,N}^{(ii)},
$$
where, using the notation introduced in (\ref{integradirichlet}), each term is defined as follows:
\begin{equation}
\label{BBi}
BB_{M,n,N}^{(i)}:=  {1\over {2\pi}}{1\over {M+1}}\sum_{|k| \leq M}\left(1-{|k|\over {M+1}}\right) k^2 \int_0^{2\pi} Y_{n,N}^2(s,s) v(s) ds,
\end{equation}
\begin{equation}
\label{BBii}
BB_{M,n,N}^{(ii)}:= - \int_0^{2\pi} {1\over {2\pi}}\int_0^u {1\over {M+1}} F^{\prime\prime}_{M,n}(u-s)Y_{n,N}(s,s) \sigma(s)dW_s Y_{n,N}(u,u) \sigma(u)dW_u,
\end{equation}

\begin{equation}
\label{BCi}
BC_{M,n,N}^{(i)}:= - \int_0^{2\pi} {1\over {2\pi}}\int_0^s {1\over {M+1}} F^{\prime \prime}_{M,n}(s-v)D_{N,n}(s-v)\sigma(v)dW_v \, Y_{n,N}(s,s) v(s)ds,
\end{equation}
\begin{equation}
\label{BCii}
BC_{M,n,N}^{(ii)}:= {1\over {M+1}}\sum_{|k| \leq M} \left(1-{|k|\over {M+1}}\right) k^2 {1\over {2\pi}}\int_0^{2\pi} \, \int_0^u e^{-ik\phi_n(s)} Y_{n,N}(s,s)\sigma(s) dW_s \, Y_{k,n,N}(u,u)\,  \sigma(u) dW_u,
\end{equation}
\begin{equation}
\label{BCiii}
BC_{M,n,N}^{(iii)}:= - \int_0^{2\pi} \int_0^u  {1\over {2\pi}}\int_0^{v} {1\over {M+1}}F^{\prime \prime}_{M,n}(v-s) D_{N,n}(v-s)\sigma(s) dW_{s} \, \sigma(v) dW_v \,
Y_{N,n}(u,u) \sigma(u) dW_u,
\end{equation}
\begin{equation}
\label{CCi}
CC_{M,n,N}^{(i)}:= {1\over {2\pi}}{1\over {M+1}}\sum_{|k| \leq M} \left(1-{|k|\over {M+1}}\right)k^2 \int_0^{2\pi} Y_{k,n,N}(s,s)Y_{-k,n,N}(s,s) v(s) ds,
\end{equation}
\begin{equation}
\label{CCii}
CC_{M,n,N}^{(ii)}:= {1\over {M+1}}\sum_{|k| \leq M} \left(1-{|k|\over {M+1}}\right) k^2  {1\over {2\pi}}\int_0^{2\pi} \int_0^s Y_{k,n,N}(v,v) \sigma(v) dW_v \, Y_{-k,n,N}(s,s) \sigma(s) dW_s.
\end{equation}

\bigskip

In summary, the estimation error
\begin{equation}
\label{maindifference}
\widehat\gamma^2_{n,N,M} - {1\over {2\pi}}\int_0^{2\pi}\gamma^2(t) dt
\end{equation}
comprises the study of four main components:
\begin{equation}
\label{continuous2}
{2\pi\over {M+1}}\sum_{|k| \leq M} \left(1-{|k|\over {M+1}}\right) \, k^2 c_{k}(v)c_{-k}(v) - {\frac{1}{2\pi}}\int_0^{2\pi} \gamma^2(t)dt
\end{equation}
\begin{equation}
\label{mainpart0}
+{2\pi\over {M+1}}\sum_{|k| \leq M} \left(1-{|k|\over {M+1}}\right) k^2 \, \left(A_{k,n}A_{-k,n}-  c_{k}(v)c_{-k}(v)\right)
\end{equation}
\begin{equation}
\label{mainpartVar}
+  2 \, (AB^{(i)}_{M,n,N}+AB^{(ii)}_{M,n,N}+AC^{(i)}_{M,n,N}+AC^{(ii)}_{M,n,N} + BB_{M,n,N}^{(ii)}+ BC_{M,n,N}^{(ii)}+BC_{M,n,N}^{(iii)}+CC_{M,n,N}^{(ii)})
\end{equation}
\begin{equation}
\label{correction}
+BB_{M,n,N}^{(i)} +2 BC_{M,n,N}^{(i)}+CC_{M,n,N}^{(i)} - K\widehat \sigma^4_{n,N,M},
\end{equation}
where $\widehat \sigma^4_{n,N,M}$ is defined in (\ref{ESTQUART}) and the constant $K$ is determined in (\ref{constantBiasCorrect}).

Accordingly, the proof of the theorem is divided into four steps. The first and second steps comprise the study of the bias correction term (\ref{correction}) and the asymptotic negligibility of the discretization error (\ref{mainpart0}). The other two steps follow \cite{JACO97} in order to identify the asymptotic variance and prove the stable convergence in law.
\par\noindent
In the proof we consider the case of regular sampling, i.e., $\varphi_n(t)={2\pi\over n}j$, \, if ${2\pi\over n}j \leq t < {2\pi \over n}(j+1)$, $j=0, \ldots, n$. Further, $C$ will always denote a positive constant, not necessarily the same.

\subsection{Step I. The bias-correction term}
\label{biascorrectionterm}

Firstly, we show that the term (\ref{correction}) is $o_p(\rho(n)^{1/4})$, therefore proving that the error (\ref{maindifference}) is equal to

$$
{2\pi\over {M+1}}\sum_{|k| \leq M} \left(1-{|k|\over {M+1}}\right) \, k^2 c_{k}(v)c_{-k}(v) -  {\frac{1}{2\pi}}\int_0^{2\pi} \gamma^2(t)dt
$$
$$
+{2\pi\over {M+1}}\sum_{|k| \leq M} \left(1-{|k|\over {M+1}}\right) k^2 \, \left(A_{k,n}A_{-k,n}-  c_{k}(v)c_{-k}(v)\right)
$$
$$
+2 \, (AB^{(i)}_{M,n,N}+AB^{(ii)}_{M,n,N}+AC^{(i)}_{M,n,N}+AC^{(ii)}_{M,n,N}
+ BB_{M,n,N}^{(ii)}+ BC_{M,n,N}^{(ii)}+BC_{M,n,N}^{(iii)}+CC_{M,n,N}^{(ii)})
$$
$$
+ o_p(\rho(n)^{1/4}).
$$
We begin by studying the term $BB_{M,n,N}^{(i)}$ defined by (\ref{BBi}).
The term $BC_{M,n,N}^{(i)}$ defined by (\ref{BCi}) and the term $CC_{M,n,N}^{(i)}$ defined by (\ref{CCi}) are analogous to $BB_{M,n,N}^{(i)}$
and give the same contribution.

Using the It\^o formula, the term
$$
BB_{M,n,N}^{(i)}:={1\over {M+1}}\sum_{|k| \leq M}\left(1-{|k|\over {M+1}}\right) k^2{1\over {2\pi}} \int_0^{2\pi} Y_{n,N}^2(s,s) v(s) ds
$$
reads as
\begin{equation}
\label{BBi1}
{1\over {M+1}}  \sum_{|k| \leq M}\left(1-{|k|\over {M+1}}\right) k^2 {1\over {2\pi}}\int_0^{2\pi} \int_0^s D_{N,n}^2(s-u) v(u) du \, v(s) ds
\end{equation}
\begin{equation}
\label{BBi2}
+{1\over {M+1}}  \sum_{|k| \leq M}\left(1-{|k|\over {M+1}}\right) k^2 {1\over {2\pi}}\int_0^{2\pi} 2 \int_0^s Y_{n,N}(u,s) D_{N,n}(s-u) \sigma(u) dW_u \, v(s) ds.
\end{equation}
The leading term is the first one, namely (\ref{BBi1}), which is easily seen to be equal to
\begin{equation}
\label{BIAScomput}
{1\over {M+1}} {1\over {6}} M(M+1)(M+2)  \, {1\over {2\pi}}\int_0^{2\pi} \, \int_0^s D_{N,n}^2(s-u) v(u) du \, v(s) ds.
\end{equation}
Now, using Lemma \ref{LemmaDirichlet1} and noting that $N/n \sim c_N/(2\pi)$ and $M^2/n \sim c^2_M/(2\pi)$,
we have that (\ref{BIAScomput}) converges in probability to
$$
{c_M^2 \over {\textcolor{black}{12}}} \, (1+2\eta({c_N/ \pi})) {1\over {2\pi}} \int_0^{2\pi} \sigma^4(t) dt.
$$

Consider now (\ref{BBi2}). Exploiting the boundedness of the process $v$, it is enough to observe that
$$
E\left[\left(\int_0^s Y_{n,N}(u,s) D_{N,n}(s-u) \sigma(u) dW_u \right)^2\right] =E\left[\int_0^s Y^2_{n,N}(u,s) D^2_{N,n}(s-u) v(u) du\right]$$
$$
\leq C \,  \int_0^s \int_0^u D^2_{N,n}(s-s_1) ds_1 \,  D^2_{N,n}(s-u)  du = O(\rho(n)^3),
$$
where we have used Lemma \ref{LemmaDirichlet1} and the property that $D^2_{N,n}(s-s_1) \leq C N^{-2}$ for $s_1< s-\varepsilon$, $\varepsilon >0$, for $n$ large enough. Therefore, in probability
the term (\ref{BBi2}) has order $M^2 \rho(n)^{3/2}\sim c_M^2 \rho(n)^{-1}\rho(n)^{3/2}= c_M^2 \rho(n)^{1/2}$, and thus it converges to zero.

\medskip

Finally, let
\begin{equation}
\label{constantBiasCorrect}
K:= \left({1\over {\textcolor{black}{12}}}+ {2\over {\textcolor{black}{12}}}+{1\over {\textcolor{black}{12}}}\right) \, {c_M^2\over {2\pi}}  (1+2\eta(c_N/\pi)) = {\textcolor{black}{1} \over {3}} \, {c_M^2\over {2\pi}} (1+2\eta(c_N/\pi)).
\end{equation}
Then, following \cite{LMM}, Theorem 3, the following convergence in probability holds:
$$
 \rho(n)^{-1/2}\left(K \hat \sigma^4_{n,N,M}- (BB_{M,n,N}^{(i)} +2 BC_{M,n,N}^{(i)}+CC_{M,n,N}^{(i)})\right) \to K( X_{c_M}+Y_{c_M,c_N}),
$$
where $X_{c_M}$ are $Y_{c_M,c_N}$ are defined as
$$X_{c_M}:= -{1\over c_M} {1\over \pi} \int_0^{2\pi} \gamma^2(t) dt, \ \ \ \ Y_{c_M,c_N}:= 2c_M {1\over \pi}(1+2\eta(c_N/\pi)) \int_0^{2\pi} \sigma^4(t) dt.$$
Therefore, the proof that (\ref{correction}) is $o_p(\rho(n)^{1/4})$ is complete.

\subsection{Step II. The statistical error}
\label{StatError}

In this paragraph we consider the discretization error component (\ref{mainpart0}), which reads as
$$
{ 1\over {M+1}} {1\over {2\pi}}\left(\int_0^{2\pi} 2 \int_0^s F_{M,n}^{\prime \prime}(s-u)v(u) du \, v(s) ds -
\int_0^{2\pi} 2 \int_0^s F_{M}^{\prime \prime}(s-u)v(u) du \, v(s) ds \right)
$$
\begin{equation}
\label{AAnew1}
={ 1\over {M+1}} {1\over {2\pi}}\int_0^{2\pi} 2 \int_0^s \left(F_{M,n}^{\prime \prime}(s-u)-F_{M}^{\prime \prime}(s-u)\right)\, v(u) du \, v(s) ds,
\end{equation}
and prove that it is $o_p(\rho(n)^{1/4})$.

Up to a negligible multiplicative constant, we rewrite the latter as follows\footnote{With abuse of notation, in (\ref{DEdec}) we denoted by $c_k(F^{\prime\prime}_{M,n})$ the $k$-th Fourier coefficient of $\widetilde F^{\prime\prime}_{M,n}$, which is the quantity defined as $\widetilde F^{\prime\prime}_{M,n}(t):= F^{\prime\prime}_{M }(\varphi_n(t))$. A straightforward
but lengthy proof, which is available from the authors, shows that the difference
between the two kernels is negligible.}:
$$
E^{(1)}_{M,n}+E^{(2)}_{M,n},
$$
where
\begin{equation}
\label{DEdec}
E^{(1)}_{M,n}:= {1\over {M+1}}\sum_{|k| \leq M} \left(c_k(F^{\prime\prime}_{M,n}) -c_k(F^{\prime\prime}_{M})\right)|c_k(v)|^2
\end{equation}
and
\begin{equation}
\label{DEdecc}
E^{(2)}_{M,n}:= {1\over {M+1}}\sum_{ |k | > M} c_k(F^{\prime\prime}_{M,n}) |c_k(v)|^2.
\end{equation}
Note that (see, also, Remark \ref{boundary}), the term $|c_k(v)|^2$ reads
\begin{equation}
\label{dec1}
|c_k(v)|^2 = \frac{1}{4 \pi^2 k^2} \left(\int_{0}^{2\pi} \gamma^2(s)\,ds + 2 \int_{0}^{2\pi} \int_{0}^{s} \cos(k(s-u))\gamma(u)\,dZ_u\,\gamma(s)\,dZ_s\right).
\end{equation}
Further, we compute the $k$-th Fourier coefficient of $F^{\prime\prime}_{M,n}$. First, note that, for any $j$, one has that
$$
c_k\left(I_{[\frac{2\pi}{n}j,\frac{2\pi}{n}(j+1)[}\right)=
\frac{1}{2\pi{\rm i}k} e^{-2\pi{\rm i}\frac{k}{n}j}\left(1-e^{-2\pi{\rm i}\frac{k}{n}}\right).
$$
Therefore, it holds that
\begin{equation}
\label{coeff2ndderfej}
c_k(F^{\prime\prime}_{M,n})
= \sum_{j=0}^{n-1}  F^{\prime\prime}_M\left(\frac{2\pi}{n}\,j\,\right) c_k\left(I_{[\frac{2\pi}{n}j,\frac{2\pi}{n}(j+1)[}\right)
\end{equation}
$$
= - \frac{1}{2\pi{\rm i} k}  \left(1-e^{-2\pi{\rm i}\frac{k}{n}} \right)\sum_{|l| \leq M} \left(1-\frac{|l|}{M+1}\right) \, l^2
\sum_{j=0}^{n-1} e^{-2\pi{\rm i} \frac{k-l}{n}j}.
$$

Now  observe that,  if $n$ divides $k,$ then   $1 - e^{-2\pi {\rm i} \frac{k}{n}}$ is equal to zero; accordingly, we assume that $k = nq +r, \, r \neq 0$  with  either  $q=0$   if $|k|\le n$ or $q \neq 0$ otherwise.  Moreover, note that the summation in (\ref{coeff2ndderfej}) is either equal to $n$, if $n$ divides $k-l,$ or equal to zero, otherwise; hence, we set $l= r $, with $|r|\le  M.$ Thus, (\ref{coeff2ndderfej}) reduces to
$$
c_k(F^{\prime\prime}_{M,n}) = - \left(1 - \frac{{|r|}}{M+1}\right)\,r^2\,\frac{1}{2 \pi {\rm i} k}\,n\,(1 - e^{- 2\pi{\rm i} \frac{r}{n}}).
$$

We now study the asymptotic behavior of the terms $E^{(1)}_{M,n}$ and $E^{(2)}_{M,n}$, separately.
First, we prove that  $E^{(2)}_{M,n}$  converges to zero in the $L^{1}$-norm. By taking into account the decomposition of $|c_k(v)|^2$ in (\ref{dec1}), we have to prove that both terms resulting from such a decomposition are asymptotically negligible. However, here
we explicitly compute the upper bound for the first term, which is the leading term.
Using the bound $\vert 1-e^{-2 \pi{\rm i} \frac{r}{n}}\vert \leq 2 \pi  \frac{|r|}{n}$, it holds that
\begin{equation}
\label{boundHF}
E\left[\left|E^{(2)}_{M,n}\right|\right]= \frac{1}{M+1}  E\left[\left| 2\sum_{k> M} c_k(F^{\prime\prime}_{M,n}) \frac{1}{4\pi k^2} \int_{0}^{2\pi} \gamma^2(s)\,ds \right| \right]
\end{equation}
$$
\le  C \frac{n}{M+1}  E\left[\int_{0}^{2\pi} \gamma^2(s)\,ds \right]\sum_{k> M} \frac{r^2}{k^3} \left(1 - \frac{{r}}{M+1}\right)\,
 \vert 1 - e^{- 2\pi {\rm i} \frac{r}{n}}  \vert
$$
$$
\le \frac{C}{M+1}   \sum_{k>M} \frac{r^3}{k^3} \left(1 - \frac{{r}}{M+1}\right)
= \frac{C}{(M+1)^2} \sum_{q = 1}^{\infty} \sum_{r = 1}^{M}  \frac{r^3(M+1 - r)}{(nq)^3(1 + \frac{r}{n q})^3}
$$
$$
\le \frac{C  {M^4}}{n^3(M+1)^2} \sum_{q = 1}^{\infty} \frac{1}{q^3}\sum_{r = 1}^{M} \left(1 + \frac{r}{nq}\right)^{-3}.
$$
Now, note that
$$
\sum_{r = 1}^{M} \left(1 + \frac{r}{n q}\right)^{-3}  \leq \frac{1}{2}n q\left(1 -  \left(1 + \frac{M}{n q}\right)^{-2} \right) =  (M+1) -  \frac{3}{2}\frac{(M+1)^2}{n q} + o(1).
$$
Therefore,  (\ref{boundHF}) is $O(n^{-1})$.

For what concerns $E^{(1)}_{M,n}$, we also study the convergence in the $L^1$-norm.  Again, we focus on the leading term from the decomposition in (\ref{dec1}).
For any fixed $k$, it holds:
$$
c_k(F^{\prime\prime}_{M,n}) -c_k(F^{\prime\prime}_{M})
=  -  \left(1-\frac{|k|}{M+1}\right)\,k^2\,\left(\frac{1-e^{-2\pi{\rm i}\frac{k}{n}}}{2\,\pi\,{\rm i}\,\frac{k}{n}}-1\right)
$$
$$
= - \left(1-\frac{|k|}{M+1}\right)\,k^2\,\left({\rm i}\,\pi\,\frac{k}{n}-\frac{2}{3}\pi^2\frac{k^2}{n^2} + O\left(\frac{k^3}{n^3}\right)\right).
$$
Therefore, we have that
$$
E[|E^{(1)}_{M,n}|]= \frac{1}{M+1}E\left[\left|\sum_{|k| \leq M} \left(c_k(F^{\prime\prime}_{M,n}) -c_k(F^{\prime\prime}_{M})\right)\frac{1}{4\pi k^2}\int_{0}^{2\pi}\gamma^2(s)\,ds \right| \right]  \\
$$
$$
\leq \frac{1}{M+1} \sum_{|k| \leq M} |c_k(F^{\prime\prime}_{M,n}) -c_k(F^{\prime\prime}_{M})| \frac{1}{4\pi^2 k^2}E\left[\int_{0}^{2\pi}\gamma^2(s)\,ds  \right]
$$
$$
\leq \frac{C}{n(M+1)} \sum_{|k| \leq M} \left(1-\frac{|k|}{M+1}\right)  |k|
= \frac{2C }{n(M+1)}  \frac{M(M-2)}{6}  =O(\rho(n)^{1/2}).
$$
This completes the proof of the asymptotical negligibility of the statistical error with faster rate than $\rho(n)^{1/4}$.

\subsection{Step III. Asymptotic variance}

This section follows \cite{JACO97} in order to identify the asymptotic variance and prove the stable convergence in law.

First, consider the term {(\ref{continuous2})}, namely:
\begin{equation}
\label{AA2a}
{2\pi\over {M+1}}\sum_{|k| \leq M} \left(1-{|k|\over {M+1}}\right) \, k^2 c_{k}(v)c_{-k}(v) - {1\over {2\pi}}\int_0^{2\pi} \gamma^2(t)dt.
\end{equation}
Using the integration by parts formula and Remark \ref{boundary}, it holds that
$c_k(v)={1\over {{\rm i}k}} c_k(dv)$ and (\ref{AA2a}) is equal to
$$
{2\pi\over {M+1}}\sum_{|k| \leq M} \left(1-{|k|\over {M+1}}\right) \, c_{k}(dv)c_{-k}(dv) - {1\over {2\pi}} \int_0^{2\pi} \gamma^2(t)dt.
$$
By applying the It\^o formula, the term (\ref{AA2a}) is equal to
$2 A_M(2\pi)$, where\footnote{To simplify the notation, in the following we will always omit the argument when it is equal to $2\pi$, so we will write $A_M$ instead of $A_M(2\pi)$. Similarly  for the processes in (\ref{VV}).}
\begin{equation}
\label{AA2c}
A_M(u):={1\over {M+1}}\, {1\over {2\pi}} \int_0^{u} \int_0^t F_M(t-s) dv(s)\, dv(t).
\end{equation}

Then, according to \cite{JACO97}, we determine the variance of the asymptotic distribution by studying
\begin{equation}
\label{VVvar}
\langle \rho(n)^{-1/4} 2\,V_{M,n,N}, \rho(n)^{-1/4} 2\, V_{M,n,N} \rangle_{2\pi} \, ,
\end{equation}
where
\begin{equation}
\label{VV}
V_{M,n,N}:=A_M+AB_{M,n,N}^{(i)}+AB_{M,n,N}^{(ii)}+AC_{M,n,N}^{(i)}+AC_{M,n,N}^{(ii)}
\end{equation}
$$
+BB_{M,n,N}^{(ii)}+BC_{M,n,N}^{(ii)}+BC_{M,n,N}^{(iii)}+CC_{M,n,N}^{(ii)}.
$$

In the first step we study the bracket 
\begin{equation}
\label{primo}
\rho(n)^{-1/2}\langle 2\, A_M, 2\, A_M \rangle_{2\pi}.
\end{equation}
$$
{M\over {c_M}} 4 \, {1\over {(2\pi)^2}}\int_0^{2\pi} \left({1\over {M+1}}   \int_0^t F_M(t-s) dv(s) \right)^2 \gamma^2(t) dt,
$$
which, by using the It\^o formula, is equal to
\begin{equation}
\label{F14}
{M\over {c_M}}\,  {1\over {\pi^2}}\int_0^{2\pi} {1\over {(M+1)^2}}\int_0^t F^2_M(t-s) \gamma^2(s) ds \, \gamma^2(t) dt
\end{equation}
\begin{equation}
\label{F15}
+ {M\over {c_M}} \, {2\over {\pi^2}}\int_0^{2\pi} {1\over {(M+1)^2}} \int_0^t \int_0^s F_M(t-u) dv(u) \, F_M(t-s) dv(s) \, \gamma^2(t) dt.
\end{equation}
Using Lemma \ref{LemmaFejer}, equation (\ref{squarefejerTEIC}), it is seen that the term (\ref{F14}) converges in probability to
\begin{equation}
\label{F14bis}
{1\over {2\pi}} \int_0^{2\pi} {4\over 3} \, {1\over {c_M}}\, \gamma^4(t) dt.
\end{equation}
Further, in order to prove that the term (\ref{F15}) is $o_p(1)$, it is enough to compute
$$
E\left[\left( {1\over {(M+1)^2}} \int_0^t \int_0^s F_M(t-u) dv(u) \, F_M(t-s) dv(s)\right)^2\right]^{1/2}
$$
$$
\leq C\, {1\over {M+1}} \left( {1\over {M+1}}\int_0^{t} E\left[ {1\over {M+1}} \int_0^s F^2_M(t-u) \gamma^2(u)du\right] \, F^2_M(t-s) ds \right)^{1/2}
$$
and apply Lemma \ref{LemmaFejer}, which gives that $({1/ {M}}) \int_0^s F^2_M(t-u) \gamma^2(u)du=o_p(1)$ and $({1/ {M}})\int_0^{t} F^2_M(t-s) ds=O(1)$,
because $(1/M) F^2_M(\cdot)$ is a {\sl good kernel}, according to Lemma \ref{LemmaFejer}.

\bigskip

The second step for identifying the asymptotic variance is to study
\begin{equation}
\label{ACvar2}
\rho(n)^{-1/2}\langle 2(AB_{M,n,N}^{(i)}+AB_{M,n,N}^{(ii)}+AC_{M,n,N}^{(i)}+AC_{M,n,N}^{(ii)}) ,
 2(AB_{M,n,N}^{(i)}+AB_{M,n,N}^{(ii)}+AC_{M,n,N}^{(i)}+AC_{M,n,N}^{(ii)}) \rangle_{2\pi}.
\end{equation}
The bracket contains $16$ terms giving the same contribution for symmetry.

We consider the term $AB_{M,n,N}^{(i)}$ defined by (\ref{ABi2}).
It holds: 
\begin{equation}
\label{ABi2a}
\langle \rho(n)^{-1/4}\, 2AB_{M,n,N}^{(i)}, \rho(n)^{-1/4}\, 2AB_{M,n,N}^{(i)}\rangle_{2\pi}
\end{equation}
$$
=\rho(n)^{-1/2}\, {{4\over {(2\pi)^2}}}\int_0^{2\pi}  \left(\int_0^s {1\over {M+1}} F_{M,n}^{\prime \prime}(s-u) v(u) du \right)^2 \, Y^2_{n,N}(s,s) v(s) ds
$$
\begin{equation}
\label{ABi2c}
=\rho(n)^{-1/2}\, {{1\over {\pi^2}}}\int_0^{2\pi}  \left(\int_0^s {1\over {M+1}} F_{M,n}^{\prime \prime}(s-u) v(u) du \right)^2 \, \int_0^s D^2_{N,n}(s-u) v(u) du \, v(s) ds
\end{equation}
\begin{equation}
\label{ABi2d}
+\rho(n)^{-1/2}\, {{2\over {\pi^2}}}\int_0^{2\pi}  \left(\int_0^s {1\over {M+1}} F_{M,n}^{\prime \prime}(s-u) v(u) du \right)^2 \, \int_0^s\int_0^u D_{N,n}(s-r)\sigma(r)dW_{r} \, D_{N,n}(s-u) \sigma(u) dW_u \, v(s) ds.
\end{equation}
Consider (\ref{ABi2c}). Thanks to Lemma \ref{LemmaDirichlet1}, equation (\ref{KEY}), and the result obtained in Step II for the term (\ref{AAnew1}), by virtue of the fact that $M/n \to 0$,  it is equivalent to study
\begin{equation}
\label{ABi2e}
\rho(n)^{-1/2}\, {{1\over {\pi^2}}} \, n^{-1 } \pi (1+2\eta(c_N/\pi)) \,\int_0^{2\pi}  \left(\int_0^s {1\over {M+1}} F_{M}^{\prime \prime}(s-u) v(u) du \right)^2 \, v^2(s) ds +o_p(1).
\end{equation}
Using the integration by parts formula and thanks to the boundedness of the volatility process $v$, one has that
\begin{equation}
\label{keyformula}
\int_0^s {1\over {M+1}} F_{M}^{\prime \prime}(s-u) v(u) du= \int_0^s {1\over {M+1}} F_{M}^{\prime}(s-u) dv(u)+ O_p(\rho(n)^{1/2}).
\end{equation}
Therefore (\ref{ABi2e}) gives
$$
\rho(n)^{-1/2}\, {{1\over {\pi^2}}}\pi \, n^{-1 } (1+2\eta(c_N/\pi)) \,\int_0^{2\pi}  \left(\int_0^s {1\over {M+1}} F_{M}^{\prime}(s-u) dv(u) \right)^2 \, v^2(s) ds +o_p(1)
$$
\begin{equation}
\label{ABi2f}
=\rho(n)^{-1/2}\, {1\over {\pi}} \, n^{-1 } (1+2\eta(c_N/\pi)) \, \int_0^{2\pi}  \int_0^s {1\over {(M+1)^2}} |F_{M}^{\prime}(s-u)|^2 \gamma^2(u) du \, v^2(s) ds
\end{equation}
\begin{equation}
\label{ABi2g}
+ \rho(n)^{-1/2}\, {{2\over {\pi}}} \, n^{-1 } (1+2\eta(c_N/\pi)) {1\over {(M+1)^2}}\int_0^{2\pi}  \int_0^s \int_0^u F_{M}^{\prime}(s-r) dv(r) \,  F_{M}^{\prime}(s-u) dv(u)
\, v^2(s) ds  +o_p(1).
\end{equation}
Consider (\ref{ABi2f}). It is enough to observe that, by Lemma \ref{LemmaFejer}, it holds:
$$
{1\over {M^3}} \int_0^{2\pi} \int_0^{2\pi}  |F_{M}^{\prime}(s-u)|^2 \gamma^2(u) du \, v^2(s) ds \to {2\over {15}} \pi \int_0^{2\pi}\gamma^2(s)   \, v^2(s) ds
$$
and
$$\rho(n)^{-1/2}\, {1\over {\pi}} \, n^{-1 } (1+2\eta(c_N/\pi)){M^3\over {(M+1)^2}} \to {1\over  {{\pi^2}}} {c_M \over 2} (1+2\eta(c_N/\pi)).$$
Then the term (\ref{ABi2f}) converges to
\begin{equation}
\label{ABi2f2}
{1\over {15}}  \, c_M \, (1+2\eta(c_N/\pi)) {1\over {2\pi}}\int_0^{2\pi}  \gamma^2(s) \, v^2(s) ds.
\end{equation}
It remains to prove that (\ref{ABi2g}) is asymptotically negligible. To this aim, consider
$$
E\left[\left(\int_0^s \int_0^u F_{M}^{\prime}(s-r) dv(r) \,  F_{M}^{\prime}(s-u) dv(u)\right)^2\right]^{1/2}
\leq C \left(\int_0^s |F_{M}^{\prime}(s-u)|^2 E\left[ \int_0^u |F_{M}^{\prime}(s-r)|^2 \gamma^2(r) dr\right]  \, du \right)^{1/2}.
$$
Finally, using the fact that ${\{K_M\}}_M$ is a family of {\sl good} kernels by Lemma \ref{LemmaFejer}, it holds that
${1\over {M^3}}\int_0^u |F_{M}^{\prime}(s-r)|^2 \gamma^2(r) dr =o_p(1)$,
for $r <s -\varepsilon$, $\varepsilon >0$, and
${1\over {M^3}}\int_0^s |F_{M}^{\prime}(s-u)|^2 \gamma^2(u) du=O_p(1)$. Therefore, it follows that (\ref{ABi2g}) has order
$n^{-1/2} \, M^{-2} \, M^3 \, o_p(1)= o_p(1)$.

The term (\ref{ABi2d}) goes to zero. This result follows from the fact that $D^2_{N,n}(s-r) < C\, N^{-2}$ for $r< s-\varepsilon$, for any $\varepsilon >0$ and $n$ large enough, and the same procedure followed for (\ref{ABi2e}).
Thus the term (\ref{ACvar2}) converges to
$$
{16\over {15}} \, c_M (1+2\eta(c_N/\pi)) {1\over {2\pi}}\int_0^{2\pi}  \gamma^2(s) \, v^2(s) ds.
$$

\medskip

The last contribution to the variance of the asymptotic distribution is obtained by studying the bracket
\begin{equation}
\label{secondoBBCC}
\langle \rho(n)^{-1/4} 2(BB_{M,n,N}^{(ii)}+BC_{M,n,N}^{(ii)}+BC_{M,n,N}^{(iii)}+CC_{M,n,N}^{(ii)}) ,
\rho(n)^{-1/4} 2(BB_{M,n,N}^{(ii)}+BC_{M,n,N}^{(ii)}+BC_{M,n,N}^{(iii)}+CC_{M,n,N}^{(ii)})  \rangle_{2\pi}.
\end{equation}
This bracket yields 16 terms, each giving the same contribution. Consider
\begin{equation}
\label{BBiiVar}
\langle \rho(n)^{-1/4} 2BB_{M,n,N}^{(ii)}, \rho(n)^{-1/4} 2BB_{M,n,N}^{(ii)} \rangle_{2\pi}
\end{equation}
$$
= \rho(n)^{-1/2} 4 \int_0^{2\pi}  \left({1\over {2\pi}} \int_0^u {1\over {M+1}} F^{\prime\prime}_{M,n}(u-s)Y_{n,N}(s,s)\sigma(s)dW_s \right)^2 Y^2_{n,N}(u,u) v(u) du.
$$
By applying the It\^o formula, we have that
$$
\left( \int_0^u {1\over {M+1}} F^{\prime\prime}_{M,n}(u-s)Y_{n,N}(s,s)\sigma(s)dW_s \right)^2= \int_0^u {1\over {(M+1)^2}} |F^{\prime\prime}_{M,n}(u-s)|^2 Y^2_{n,N}(s,s)v(s)ds
$$
$$
+ 2\int_0^u \int_0^s {1\over {M+1}} F^{\prime\prime}_{M,n}(u-s_1)Y_{n,N}(s_1,s)\sigma(s_1)dW_{s_1} \, {1\over {M+1}} F^{\prime\prime}_{M,n}(u-s)Y_{n,N}(s,s)\sigma(s)dW_{s}.
$$
Therefore, the bracket (\ref{BBiiVar}) splits as
\begin{equation}
\label{BBiiVar1}
\rho(n)^{-1/2} \int_0^{2\pi} {4\over {(2\pi)^2}}\int_0^u {1\over {(M+1)^2}} |F^{\prime\prime}_{M,n}(u-s)|^2Y^2_{n,N}(s,s)v(s)ds \, Y^2_{n,N}(u,u) v(u) du
\end{equation}
\begin{equation}
\label{BBiiVar2}
+\rho(n)^{-1/2}  \int_0^{2\pi} {8\over {(2\pi)^2}}\int_0^u \int_0^s {1\over {(M+1)^2}} F^{\prime\prime}_{M,n}(u-r)Y_{n,N}(r,s)\sigma(r)dW_{r} \, F^{\prime\prime}_{M,n}(u-s)Y_{n,N}(s,s)\sigma(s)dW_{s}
\end{equation}
$$\times \, Y^2_{n,N}(u,u) v(u) du.$$
Consider first (\ref{BBiiVar1}), which contributes to the asymptotic variance. As we will show later,  (\ref{BBiiVar2}) goes instead to zero. By applying the It\^o formula twice, we have that
$$
\rho(n)^{-1/2} \int_0^{2\pi} {1\over {\pi^2}}\int_0^u {1\over {(M+1)^2}} |F^{\prime\prime}_{M,n}(u-s)|^2Y^2_{n,N}(s,s)v(s)ds \, Y^2_{n,N}(u,u) v(u) du
$$
\begin{equation}
\label{BBiiVar1bis}
=\rho(n)^{-1/2}  \int_0^{2\pi} {1\over {\pi^2}}\int_0^u {1\over {(M+1)^2}} |F^{\prime\prime}_{M,n}(u-s)|^2 \int_0^s D^2_{N,n}(s-u^{\prime}) v(u^{\prime}) du^{\prime} \, v(s)ds \,  \int_0^u D^2_{N,n}(u-s^{\prime}) v(s^{\prime}) ds^{\prime} \, v(u) du
\end{equation}
\begin{equation}
\label{BBiiVar1ter}
+\rho(n)^{-1/2}  \int_0^{2\pi} {2\over {\pi^2}}\int_0^u {1\over {(M+1)^2}} |F^{\prime\prime}_{M,n}(u-s)|^2 \int_0^s\int_0^{s^{\prime}} D_{N,n}(s-r)\sigma(r)dW_r \, D_{N,n}(s-s^{\prime}) \sigma(s^{\prime})dW_{s^{\prime}}\,  v(s)ds
\end{equation}
$$
\times \int_0^u D^2_{N,n}(u-s^{\prime}) v(s^{\prime}) ds^{\prime} \, v(u) du
$$
\begin{equation}
\label{BBiiVar1quater}
+\rho(n)^{-1/2}  \int_0^{2\pi} {2\over {\pi^2}}\int_0^u {1\over {(M+1)^2}} |F^{\prime\prime}_{M,n}(u-s)|^2 \int_0^s D^2_{N,n}(s-u^{\prime}) v(u^{\prime}) du^{\prime} \, v(s)ds
\end{equation}
$$
\times  \, \int_0^u\int_0^{u^{\prime}} D_{N,n}(u-r)\sigma(r)dW_r \, D_{N,n}(u-u^{\prime}) \sigma(u^{\prime})dW_{u^{\prime}}\,  v(u)du
$$
\begin{equation}
\label{BBiiVar1quinto}
+\rho(n)^{-1/2} \int_0^{2\pi} {4\over {\pi^2}}\int_0^u {1\over {(M+1)^2}} |F^{\prime\prime}_{M,n}(u-s)|^2 \int_0^s\int_0^{s^{\prime}} D_{N,n}(s-r)\sigma(r)dW_r \, D_{N,n}(s-s^{\prime}) \sigma(s^{\prime})dW_{s^{\prime}}\,  v(s)ds
\end{equation}
$$
\times  \, \int_0^u\int_0^{u^{\prime}} D_{N,n}(u-r)\sigma(r)dW_r \, D_{N,n}(u-u^{\prime}) \sigma(u^{\prime})dW_{u^{\prime}}\,  v(u)du .
$$
Using Lemma \ref{LemmaDirichlet1}, equation (\ref{KEY}), the term (\ref{BBiiVar1bis}) gives
$$
\rho(n)^{-1/2} {1\over {n^2}} (1+2\eta(c_N/\pi))^2 \, \int_0^{2\pi} \int_0^u {1\over {(M+1)^2}} |F^{\prime\prime}_{M,n}(u-s)|^2 \sigma^4(s)
ds \, \sigma^4(u) \, du.
$$
Finally, using Lemma \ref{LemmaFejer}, the term (\ref{BBiiVar1bis}) converges to
$${1\over {105}} \, c^3_M \, (1+2\eta(c_N/\pi))^2  \, {1\over {2\pi}}\int_0^{2\pi}\sigma^8(t)dt.$$
A similar procedure as for (\ref{ABi2d}) allows us to prove that (\ref{BBiiVar1ter}), (\ref{BBiiVar1quater}) and (\ref{BBiiVar1quinto}) go to zero in probability.

\bigskip

We verify now that the (\ref{BBiiVar2}) goes to zero in probability. It is enough to show that
\begin{equation}
\label{BBiiVar2new}
\int_0^u \int_0^s {1\over {(M+1)^2}}F_{M,n}^{\prime\prime}(u-r)Y_{n,N}(r,s) \sigma(r)dW_{r} \, F_{M,n}^{\prime\prime}(u-s)Y_{n,N}(s,s) \sigma(s)dW_{s} \, Y^2_{n,N}(u,u)
\end{equation}
is $o_p(\rho(n)^{1/2})$. By the It\^o formula and Lemma \ref{LemmaDirichlet1}, it holds that
$$
\int_0^u \int_0^s {1\over {(M+1)^2}}F_{M,n}^{\prime\prime}(u-r)Y_{n,N}(r,s) \sigma(r)dW_{r} \, F_{M,n}^{\prime\prime}(u-s)Y_{n,N}(s,s) \sigma(s)dW_{s} \, Y^2_{n,N}(u,u)
$$
\begin{equation}
\label{BBiiVar2a}
= \int_0^u \int_0^s {1\over {(M+1)^2}}F_{M,n}^{\prime\prime}(u-r)Y_{n,N}(r,s) \sigma(r)dW_{r} \, F_{M,n}^{\prime\prime}(u-s)Y_{n,N}(s,s) \sigma(s)dW_{s} \int_0^u D^2_{N,n}(u-r)v(r) dr
\end{equation}
$$+ o_p(\rho(n)^{1/2}).
$$
Consider (\ref{BBiiVar2a}). By the Cauchy-Schwarz inequality, we have that
$$
E\left[\left|\int_0^u \int_0^s {1\over {(M+1)^2}}F_{M,n}^{\prime\prime}(u-r)Y_{n,N}(r,s) \sigma(r)dW_{r} \, F_{M,n}^{\prime\prime}(u-s)Y_{n,N}(s,s) \sigma(s)dW_{s} \int_0^u D^2_{N,n}(u-r)v(r) dr\right|\right]
$$
$$
\leq E\left[\left(\int_0^u \int_0^s {1\over {(M+1)^2}}F_{M,n}^{\prime\prime}(u-r)Y_{n,N}(r,s) \sigma(r)dW_{r} \, F_{M,n}^{\prime\prime}(u-s)Y_{n,N}(s,s) \sigma(s)dW_{s}\right)^2\right]^{1/2}
$$
$$\times E\left[\left(\int_0^u D^2_{N,n}(u-r)v(r) dr\right)^2\right]^{1/2}.$$
Note that $\int_0^u D^2_{N,n}(u-r)v(r) dr=O_p(\rho(n))$, by Lemma \ref{LemmaDirichlet1}.
Consider then
$$
E\left[\left(\int_0^u \int_0^s {1\over {(M+1)^2}}F_{M,n}^{\prime\prime}(u-r)Y_{n,N}(r,s) \sigma(r)dW_{r} \, F_{M,n}^{\prime\prime}(u-s)Y_{n,N}(s,s) \sigma(s)dW_{s}\right)^2\right]
$$
$$
=E\left[\int_0^u \left(\int_0^s {1\over {(M+1)^2}}F_{M,n}^{\prime\prime}(u-r)Y_{n,N}(r,s) \sigma(r)dW_{r}\right)^2 \, |F_{M,n}^{\prime\prime}(u-s)|^2Y^2_{n,N}(s,s) v(s)d{s}\right]
$$
$$
=E\left[\int_0^u \int_0^s {1\over {(M+1)^4}}|F_{M,n}^{\prime\prime}(u-r)|^2 Y^2_{n,N}(r,s) v(r)d{r} \, |F_{M,n}^{\prime\prime}(u-s)|^2 Y^2_{n,N}(s,s) v(s)d{s}\right] + o(\rho(n))
$$
\begin{equation}
\label{BBiiVar2c}
=E\left[\int_0^u \int_0^s {1\over {(M+1)^4}}|F_{M,n}^{\prime\prime}(u-r)|^2\int_0^{r} D^2_{N,n}(s-u^{\prime}) v(u^{\prime})du^{\prime} \,v(r)d{r} \, |F_{M,n}^{\prime\prime}(u-s)|^2Y^2_{n,N}(s,s) v(s)d{s}\right] + o(\rho(n)).
\end{equation}
Again, using that $D^2_{N,n}(s-u^{\prime})\leq C/N^2$ for $u^{\prime} < s-\varepsilon$, for $\varepsilon >0$ and $n$ large enough, then (\ref{BBiiVar2c}) is smaller than
$$
{C\over {N^2}} {1\over {(M+1)^4}}E\left[\int_0^u |F_{M,n}^{\prime\prime}(u-s)|^2 \int_0^s |F_{M,n}^{\prime\prime}(u-r)|^2 \,v(r)d{r} \,Y^2_{n,N}(s,s)ds\right]
$$
$$
={C\over {N^2}} {1\over {(M+1)^4}}E\left[\int_0^u |F_{M,n}^{\prime\prime}(u-s)|^2 \int_0^s |F_{M,n}^{\prime\prime}(u-r)|^2 \,v(r)d{r} \,
\left(\int_0^{s} D^2_{N,n}(s-u^{\prime}) v(u^{\prime})du^{\prime}+o_p(\rho(n))\right) ds\right]
$$
$$
={C\over {N^3}} {1\over {(M+1)^4}}E\left[\int_0^u |F_{M,n}^{\prime\prime}(u-s)|^2 \int_0^s |F_{M,n}^{\prime\prime}(u-r)|^2 \,v(r)d{r} \, v(s) ds+o_p(\rho(n))\right].
$$
Then, using Lemma \ref{LemmaFejer}, and the fact that $({1/ {M^5}}) |F_{M,n}^{\prime\prime}(x)|^2$ is a {\sl good} kernel, then this term has order
$N^{-3}M^6 o(1)= o(1)$.
Finally, we obtain the order of (\ref{BBiiVar2new}), which is
$o_p(\rho(n)^{1/2})$.

\medskip

We can conclude that the total contribution of all the terms in (\ref{secondoBBCC}) is
$${16 \over {105}} \, c^3_M \, (1+2\eta(c_N/\pi))^2  {1\over {2\pi}}\int_0^{2\pi}\sigma^8(t)dt.$$

\bigskip

It remains to show that the other brackets in (\ref{VVvar}) give asymptotically negligible contributions. We study in detail the convergence in probability
\begin{equation}
\label{AABi}
\langle \rho(n)^{-1/4} A_M, \rho(n)^{-1/4} AB^{(i)}_{M,n,N} \rangle_{2\pi} \to 0
\end{equation}
and the convergence in probability
\begin{equation}
\label{ABiBBii}
\langle \rho(n)^{-1/4} AB^{(i)}_{M,n,N}, \rho(n)^{-1/4} BB^{(ii)}_{M,n,N} \rangle_{2\pi} \to 0.
\end{equation}
The proof is analogous for the other terms. The bracket (\ref{AABi}) is equal to
\begin{equation}
\label{AABi1}
\rho(n)^{-1/2}{1\over {(2\pi)^2}} {1\over {(M+1)^2}}\int_0^{2\pi} \int_0^t F_M(t-s)dv(s)\int_0^t F_{M,n}^{\prime\prime}(t-s)v(s) ds \, Y_{n,N}(t,t) \sigma(t) \gamma(t) \rho\, dt.
\end{equation}
Omitting the negligible constants, and using the result obtained in Step II for the term (\ref{AAnew1}), by virtue of the fact that $M/n \to 0$, then (\ref{AABi1}) is equal to
$$
\rho(n)^{-1/2} \left( {1\over {(M+1)^2}}\int_0^{2\pi} \int_0^t F_M(t-s)dv(s)\int_0^t F_{M}^{\prime\prime}(t-s)v(s) ds \, Y_{n,N}(t,t) \sigma(t) \gamma(t) \rho\, dt + o_p(\rho(n)^{1/2}) \right).
$$
Moreover, by (\ref{keyformula}), we are lead to study
\begin{equation}
\label{AABi3}
{1\over {M+1}}\int_0^{2\pi} \int_0^t F_M(t-s)dv(s)\int_0^t F_{M}^{\prime}(t-s)dv(s) \, Y_{n,N}(t,t) \sigma(t) \gamma(t) \rho\, dt.
\end{equation}
By applying the It\^o formula, (\ref{AABi3}) is equal to
\begin{equation}
\label{AABi4}
{1\over {M+1}}\int_0^{2\pi} \int_0^t F_M(t-s) F_{M}^{\prime}(t-s)\gamma^2(s)ds \, Y_{n,N}(t,t) \sigma(t) \gamma(t) \rho\, dt
\end{equation}
\begin{equation}
\label{AABi5}
+ {1\over {M+1}}\int_0^{2\pi} \int_0^t \int_0^s F_M(t-u)dv(u)\, F_{M}^{\prime}(t-s)dv(s) \, Y_{n,N}(t,t) \sigma(t) \gamma(t) \rho\, dt
\end{equation}
\begin{equation}
\label{AABi6}
+{1\over {M+1}}\int_0^{2\pi} \int_0^t \int_0^s F_{M}^{\prime}(t-u)dv(u) \, F_M(t-s)dv(s) \, Y_{n,N}(t,t) \sigma(t) \gamma(t) \rho\, dt.
\end{equation}
Consider (\ref{AABi4}). It is equal to the sum
\begin{equation}
\label{AABi4a}
{1\over {M+1}}\int_0^{2\pi} \int_0^t \int_0^s D_{N,n}(t-u)\sigma(u) dW_u \, F_M(t-s) F_{M}^{\prime}(t-s)\gamma^2(s)ds \, \sigma(t) \gamma(t) \rho\, dt
\end{equation}
\begin{equation}
\label{AABi4b}
+{1\over {M+1}}\int_0^{2\pi} \int_0^t \int_0^s F_M(t-u) F_{M}^{\prime}(t-u)\gamma^2(u)du \, D_{N,n}(t-s)\sigma(s) dW_s \,\sigma(t) \gamma(t) \rho\, dt.
\end{equation}
We study (\ref{AABi4a}). The term (\ref{AABi4b}) is analogous. By the boundedness of the volatility and the volatility of volatility processes, it is enough to prove that
$$E\left[\left|{1\over {M+1}} \int_0^t \int_0^s D_{N,n}(t-u)\sigma(u) dW_u \, F_M(t-s) F_{M}^{\prime}(t-s)\gamma^2(s)ds \right|\right] \to 0.$$
The term (\ref{AABi4a}) is smaller than
$$
{C\over {M+1}} \int_0^t E\left[\left(\int_0^s D_{N,n}(t-u)\sigma(u) dW_u\right)^2\right]^{1/2} \, |F_M(t-s) F_{M}^{\prime}(t-s)| \, ds
$$
$$
\leq {C\over {M+1}} \int_0^t E\left[\int_0^s D^2_{N,n}(t-u)v(u) du\right]^{1/2} \, |F_M(t-s) F_{M}^{\prime}(t-s)| \, ds.
$$
Noting that $D^2_{N,n}(t-u) < {C/N^2}$ for $u< t-\varepsilon$, for $\varepsilon >0$, and $n$ large enough, the previous term is smaller than
$$
{1\over {M+1}} {C\over N }\left(\int_0^t |F_M(t-s)|^2 ds\right)^{1/2} \left(\int_0^t |F_{M}^{\prime}(t-s)|^2 ds\right)^{1/2}.
$$
Using Lemmas \ref{LemmaFejer}, this last term is $C \, M^{-1}n^{-1}M^{1/2}M^{3/2}= O(n^{-1/2})$.

Consider now (\ref{AABi5}). The term (\ref{AABi6}) is analogous.
It is enough to show that
$$
{1\over {M+1}}  E\left[\left|\int_0^t \int_0^s F_M(t-u)dv(u)\, F_{M}^{\prime}(t-s)dv(s) \, Y_{n,N}(t,t)\right|\right] \to 0.
$$
Therefore, using the Cauchy-Schwarz inequality, we consider
$$
{1\over {M+1}}  E\left[\left(\int_0^t \int_0^s F_M(t-u)dv(u)\, F_{M}^{\prime}(t-s)dv_s )^2]^{1/2}\, E[(Y_{n,N}(t,t)\right)^2\right]^{1/2}.
$$
Using the fact that $E[(Y_{n,N}(t,t))^2]^{1/2}=O(n^{-1/2})$, the It\^o isometry and the boundedness of volatility of volatility, we have:
$$
E\left[\left(\int_0^t \int_0^s F_M(t-u)dv(u)\, F_{M}^{\prime}(t-s)dv(s) \right)^2\right]= E\left[\int_0^t \left(\int_0^s F_M(t-u)dv(u)\right)^2\, |F_{M}^{\prime}(t-s)|^2\gamma^2(s)ds \right]
$$
$$
\leq C \int_0^t |F_{M}^{\prime}(t-s)|^2 \, E\left[\left(\int_0^s F_M(t-u)dv(u)\right)^2\right]\,ds= C \, M^4 \int_0^t {1\over {M^3}}|F_{M}^{\prime}(t-s)|^2 \, E\left[{1\over M}\int_0^s F^2_M(t-u)\gamma^2(u)du \right]\,ds.
$$
By Lemma \ref{LemmaFejer}, it holds that, for $u<t-\varepsilon$,
${1\over M}\int_0^s F^2_M(t-u)\gamma^2(u)du =o_p(1)$, and ${1\over {M^3}}\int_0^t |F_{M}^{\prime}(t-s)|^2 ds =O(1)$. Then
(\ref{AABi5}) has order $M^{-1} n^{-1/2} M^{3/2} M^{1/2} o_p(1)= o_p(1)$.

Consider now the bracket (\ref{ABiBBii}). Neglecting the irrelevant constants, this is equal to
$$
\rho(n)^{-1/2} \int_0^{2\pi} {1\over {(M+1)^2}}\int_0^s F_{M,n}^{\prime\prime}(s-u)v(u)du\int_0^s F_{M,n}^{\prime\prime}(s-u)Y_{n,N}(u,u)\sigma(u)dW_u \, Y^2_{n,N}(s,s) v(s) ds.
$$
Applying the It\^o formula, we have to study:
\begin{equation}
\label{ABiBBii1}
{1\over {M}}\int_0^{2\pi}\int_0^s \int_0^u F_{M,n}^{\prime\prime}(u-u_1) Y_{n,N}(u_1,u_1)\sigma(u_1)dW_{u_1}\, F_{M,n}^{\prime\prime}(s-u) v(u) du \, Y^2_{n,N}(s,s) v(s) ds
\end{equation}
$$
+{1\over {M}}\int_0^{2\pi}\int_0^s \int_0^u F_{M,n}^{\prime\prime}(u-u_1)v(u_1)du_1\, F_{M,n}^{\prime\prime}(s-u)Y_{n,N}(u,u)\sigma(u)dW_u \, Y^2_{n,N}(s,s) v(s) ds.
$$
We study (\ref{ABiBBii1}). The second addend is analogous. By using the boundedness of the volatility process and applying the Cauchy-Schwarz inequality, it is enough to study
\begin{equation}
\label{ABiBBii1b}
E\left[\left({1\over {M}}\int_0^s \int_0^u F_{M,n}^{\prime\prime}(u-u_1)Y_{n,N}(u_1,u_1)\sigma(u_1)dW_{u_1}\, F_{M,n}^{\prime\prime}(s-u)v(u)du\right)^2\right]^{1/2}   E[Y^4_{n,N}(s,s)]^{1/2}.
\end{equation}
Using the Burkholder-Davis-Gundy inequality and Lemma \ref{LemmaDirichlet1}, it holds that
\begin{equation}
\label{quarta}
E[Y^4_{n,N}(s,s)]^{1/2}\leq C \rho(n).
\end{equation}
Then we compute
$$
E\left[\left({1\over {M}}\int_0^s \int_0^u F_{M,n}^{\prime\prime}(u-u_1)Y_{n,N}(u_1,u_1)\sigma(u_1)dW_{u_1}\, F_{M,n}^{\prime\prime}(s-u)v(u)du\right)^2\right]^{1/2}
$$
$$
=E\left[{1\over {M^2}}\int_0^s F_{M,n}^{\prime\prime}(s-u)v(u)du \int_0^s F_{M,n}^{\prime\prime}(s-u^{\prime})v(u^{\prime})du^{\prime}  \right.
$$
$$
\left. \times \int_0^u F_{M,n}^{\prime\prime}(u-u_1)Y_{n,N}(u_1,u_1)\sigma(u_1)dW_{u_1}\,\int_0^{u^{\prime}} F_{M,n}^{\prime\prime}(u^{\prime}-u_1)Y_{n,N}(u_1,u_1)\sigma(u_1)dW_{u_1}\right]^{1/2}.
$$
For symmetry, we can assume that $u^{\prime}\leq u$ and study the following:
$$
E\left[{1\over {M^2}}\int_0^s F_{M,n}^{\prime\prime}(s-u)v(u)du \int_0^s F_{M,n}^{\prime\prime}(s-u^{\prime})v(u^{\prime})du^{\prime}
\int_0^u 1_{[0,u^{\prime}]}(u_1)F_{M,n}^{\prime\prime}(u-u_1)F_{M,n}^{\prime\prime}(u^{\prime}-u_1)
Y^2_{n,N}(u_1,u_1)v(u_1)du_1\right]^{1/2}.
$$
Moreover, using Lemma \ref{LemmaDirichlet1}, it is enough to study
\begin{equation}
\label{ABiBBii1c}
\rho(n) E\left[{1\over {M^2}}\int_0^s F_{M,n}^{\prime\prime}(s-u)v(u)du \int_0^s F_{M,n}^{\prime\prime}(s-u^{\prime})v(u^{\prime})du^{\prime}
\int_0^u 1_{[0,u^{\prime}]}(u_1)F_{M,n}^{\prime\prime}(u-u_1)F_{M,n}^{\prime\prime}(u^{\prime}-u_1) v(u_1)du_1\right].
\end{equation}
The order of (\ref{ABiBBii1c}) is studied by splitting the term as follows. For any $\varepsilon >0$, consider first the case $u^{\prime}<u-\varepsilon$ and note that
$$
E\left[{1\over {M^2}}\int_0^u 1_{[0,u^{\prime}]}(u_1)F_{M,n}^{\prime\prime}(u-u_1)F_{M,n}^{\prime\prime}(u^{\prime}-u_1) v(u_1)du_1\right] \sim M^3 o(1).
$$
Therefore, by using (\ref{keyformula}), we have obtained that (\ref{ABiBBii1b}) has order $n^{-1}n^{-1/2}M^{3/2}M^{3/2}o(1)=o(1)$.

Finally, consider the case $u-\varepsilon\leq u^{\prime}\leq u$ and, with similar arguments as above, observe that
\small
$$
\rho(n) E\left[{1\over {M^2}}\int_0^s F_{M,n}^{\prime\prime}(s-u)v(u)du \int_0^s 1_{[u-\varepsilon \leq u^{\prime}\leq u]}(u^{\prime}) F_{M,n}^{\prime\prime}(s-u^{\prime})v(u^{\prime})du^{\prime}
\int_0^u 1_{[0,u^{\prime}]}(u_1)F_{M,n}^{\prime\prime}(u-u_1)F_{M,n}^{\prime\prime}(u^{\prime}-u_1) v(u_1)du_1\right]
$$
\normalsize
is $\varepsilon \, O(n^{-1}M^{3})$. Thus, in this case, the term (\ref{ABiBBii1b}) is $\varepsilon \, O(n^{-1}n^{-1/2}M^{3/2}M^{3/2})= \varepsilon \, O(1)$, for any $\varepsilon >0$.

\subsection{Step IV. Orthogonality}
\label{Orto}

The final step requires to prove that in probability, as $n,N,M \to \infty$,
$$
\langle \rho(n)^{-1/4} 2 \, V_{M,n,N},W  \rangle_{2\pi} \to 0.
$$
We provide a detailed proof for the convergence to $0$ in probability of the bracket
\begin{equation}
\label{orthoBB}
\langle \rho(n)^{-1/4} \, BB_{M,n,N}^{(ii)} , W  \rangle_{2\pi}.
\end{equation}
The convergence of the other terms can be shown with an analogous procedure.

Consider (\ref{orthoBB}), omitting the negligible constants.
We show that
\begin{equation}
\label{ortho2}
E\left[\left(\rho(n)^{-1/4} \, \int_0^{2\pi} \int_0^u {1\over {M+1}} F^{\prime\prime}_{M,n}(u-s) Y_{n,N}(s,s) \sigma(s) dW_s \, Y_{n,N}(u,u) \sigma(u) du\right)^2\right] \to 0.
\end{equation}
Let
$$
Z_{N,M,n}(u,u):=\int_0^u {1\over {M+1}} F^{\prime\prime}_{M,n}(u-s) Y_{n,N}(s,s) \sigma(s) dW_s.
$$
Then (\ref{ortho2}) is equal to
\begin{equation}
\label{ortho3}
\rho(n)^{-1/2} E\left[ \int_0^{2\pi} \, \sigma(u)  \int_0^{2\pi} \, \sigma(u^{\prime}) Z_{N,M,n}(u,u)Z_{N,M,n}(u^{\prime},u^{\prime}) \, Y_{n,N}(u,u) Y_{n,N}(u^{\prime},u^{\prime}) \, du^{\prime} \, du\right].
\end{equation}
By symmetry, we can assume $u^{\prime}\leq u$. By the It\^o formula:
$$
Z_{N,M,n}(u,u)Z_{N,M,n}(u^{\prime},u^{\prime})
$$
\begin{equation}
\label{ortho4}
= \int_0^u 1_{[0,u^{\prime}]}(s) {1\over {(M+1)^2}} F^{\prime\prime}_{M,n}(u-s)F^{\prime\prime}_{M,n}(u^{\prime}-s)Y^2_{n,N}(s,s) v(s) ds
\end{equation}
\begin{equation}
\label{ortho4mix}
+2 \int_0^u 1_{[0,u^{\prime}]}(s) {1\over {(M+1)^2}} \int_0^s F^{\prime\prime}_{M,n}(u^{\prime}-r)Y_{N,n}(r,r) \sigma(r)dW_{r} \, F^{\prime\prime}_{M,n}(u-s)Y_{n,N}(s,s) \sigma(s) dW_s
\end{equation}
and
$$
Y_{n,N}(u,u) Y_{n,N}(u^{\prime},u^{\prime})
$$
\begin{equation}
\label{ortho5}
= \int_0^u 1_{[0,u^{\prime}]}(t) D_{N,n}(u-t) D_{N,n}(u^{\prime}-t)v(t) dt
\end{equation}
\begin{equation}
\label{ortho5mix}
+2\int_0^u 1_{[0,u^{\prime}]}(t)\int_0^t D_{N,n}(u^{\prime}-r)\sigma(r)dW_{r} \, D_{N,n}(u-t) \sigma(t)dW_t.
\end{equation}
It is enough to consider the terms in (\ref{ortho4}) and (\ref{ortho5}). In fact, the terms (\ref{ortho4mix}) and (\ref{ortho5mix}) are of higher infinitesimal order.
Substitute the terms (\ref{ortho4}) and (\ref{ortho5}) into (\ref{ortho3}) to obtain
\begin{equation}
\label{ortho6}
\rho(n)^{-1/2} E[ \int_0^{2\pi} \sigma(u)  \int_0^{2\pi} \sigma(u^{\prime}) \int_0^u 1_{[0,u^{\prime}]}(s) {1\over {(M+1)^2}} F^{\prime\prime}_{M,n}(u-s)F^{\prime\prime}_{M,n}(u^{\prime}-s)Y^2_{n,N}(s,s) v(s) ds
\end{equation}
$$
\times \int_0^u 1_{[0,u^{\prime}]}(t) D_{N,n}(u-t) D_{N,n}(u^{\prime}-t)v(t) dt \, du^{\prime} \, du].
$$
Consider
$$Y^2_{n,N}(s,s)= \int_0^s D^2_{N,n}(s-s_1) v(s_1) ds_1+ 2 \int_0^s\int_0^{s_1} D_{N,n}(s-s_2) \sigma(s_2) dW_{s_2} \, D_{N,n}(s-s_1)\sigma(s_1) dW_{s_1}$$
and plug it into equation (\ref{ortho6}). Finally, we have that
$$
\rho(n)^{-1/2} E[ \int_0^{2\pi} \sigma(u)  \int_0^{2\pi} \sigma(u^{\prime}) \int_0^u 1_{[0,u^{\prime}]}(s) {1\over {(M+1)^2}} F^{\prime\prime}_{M,n}(u-s)F^{\prime\prime}_{M,n}(u^{\prime}-s)
$$
$$
\times \int_0^s D^2_{N,n}(s-s_1) v(s_1) ds_1\, v(s) ds \, \int_0^u 1_{[0,u^{\prime}]}(t) D_{N,n}(u-t) D_{N,n}(u^{\prime}-t)v(t) dt \, du^{\prime} \, du] +o(1)
$$
\begin{equation}
\label{ortho8}
\leq C \rho(n)^{1/2}\, \int_0^{2\pi} \int_0^{2\pi}  E\left[\left| \int_0^u 1_{[0,u^{\prime}]}(s) {1\over {(M+1)^2}} F^{\prime\prime}_{M,n}(u-s)F^{\prime\prime}_{M,n}(u^{\prime}-s) v^2(s) ds \right. \right.
\end{equation}
$$
\left. \left. \times \, \int_0^u 1_{[0,u^{\prime}]}(t) D_{N,n}(u-t) D_{N,n}(u^{\prime}-t)v(t) dt\right|\right] \, du^{\prime} \, du   +o(1),
$$
where we have used the boundedness of volatility and Lemma \ref{LemmaDirichlet1}.
Then, by applying the Cauchy-Schwarz inequality, one sees that (\ref{ortho8}) is smaller than
\begin{equation}
\label{ortho9}
C \rho(n)^{1/2} \, \int_0^{2\pi}  \int_0^{2\pi}  E\left[\left( \int_0^u 1_{[0,u^{\prime}]}(s) {1\over {(M+1)^2}} F^{\prime\prime}_{M,n}(u-s)F^{\prime\prime}_{M,n}(u^{\prime}-s) v(s) ds\right)^2\right]^{1/2}
\end{equation}
$$
\times \, E\left[\left(\int_0^u 1_{[0,u^{\prime}]}(t) D_{N,n}(u-t) D_{N,n}(u^{\prime}-t)v(t) dt\right)^2\right]^{1/2} \, du^{\prime} \, du.
$$
First, consider the case $u^{\prime} < u-\varepsilon$, for any $\varepsilon >0$:
\begin{equation}
\label{ortho10}
E\left[\left( \int_0^u 1_{[0,u^{\prime}]}(s) {1\over {(M+1)^2}} F^{\prime\prime}_{M,n}(u-s)F^{\prime\prime}_{M,n}(u^{\prime}-s) v(s) ds\right)^2\right]
\end{equation}
\begin{equation}
\label{ortho11}
\leq E\left[\int_0^u 1_{[0,u^{\prime}]}(s) {1\over {(M+1)^2}} |F^{\prime\prime}_{M,n}(u-s)|^2  v(s) ds\right]
\, E\left[ \int_0^u 1_{[0,u^{\prime}]}(s) {1\over {(M+1)^2}} |F^{\prime\prime}_{M,n}(u^{\prime}-s)|^2 v(s) ds\right].
\end{equation}
Using Lemma \ref{LemmaFejer} and the fact that $u^{\prime} < u-\varepsilon$, this term has order
$M^3 \, M^3 \, o(1)$.
Then consider
\begin{equation}
\label{ortho12}
E\left[\left(\int_0^u 1_{[0,u^{\prime}]}(t) D_{N,n}(u-t) D_{N,n}(u^{\prime}-t)v(t) dt\right)^2\right]
\end{equation}
$$\leq
E\left[\int_0^u 1_{[0,u^{\prime}]}(t) D^2_{N,n}(u-t) v(t) dt\right] \,  E\left[\int_0^u 1_{[0,u^{\prime}]}(t) D^2_{N,n}(u^{\prime}-t) v(t) dt\right].
$$
Using Lemma \ref{LemmaDirichlet1} and the fact that $D^2_{N,n}(u-t)\leq {C\over {N^2}}$ for $t < u-\varepsilon$ and $n$ large enough, we see that this term has order $N^{-2} \, \rho(n)$.
Going back to (\ref{ortho9}), it has then order
$ \rho(n)^{1/2} M^{3/2} \, M^{3/2} \, N^{-1} \, \rho(n)^{1/2} \, o(1)= o(\rho(n)^{1/2}) \to 0.$

Now consider the case $u-\varepsilon \leq u^{\prime} \leq u$, for any $\varepsilon >0$. Starting from equation (\ref{ortho9}), we have
\begin{equation}
\label{ortho13}
\rho(n)^{1/2} \, \int_0^{2\pi} du  \int_0^{2\pi} 1_{[u-\varepsilon,u]}(u^{\prime}) \, du^{\prime} E\left[\left( \int_0^u 1_{[0,u^{\prime}]}(s) {1\over {(M+1)^2}} F^{\prime\prime}_{M,n}(u-s)F^{\prime\prime}_{M,n}(u^{\prime}-s) v(s) ds\right)^2\right]^{1/2}
\end{equation}
$$
\times \, E\left[\left(\int_0^u 1_{[0,u^{\prime}]}(t) D_{N,n}(u-t) D_{N,n}(u^{\prime}-t)v(t) dt\right)^2\right]^{1/2}.
$$
In this case, it is easily seen that the order of
\small
$$\rho(n)^{1/2} \, E\left[\left( \int_0^u 1_{[0,u^{\prime}]}(s) {1\over {(M+1)^2}} F^{\prime\prime}_{M,n}(u-s)F^{\prime\prime}_{M,n}(u^{\prime}-s) v(s) ds\right)^2\right]^{1/2}
  E\left[\left(\int_0^u 1_{[0,u^{\prime}]}(t) D_{N,n}(u-t) D_{N,n}(u^{\prime}-t)v(t) dt\right)^2\right]^{1/2} $$
is
$\rho(n)^{1/2} M^{3/2}M^{3/2}\rho(n)^{1/2}\rho(n)^{1/2}=O(1).$
\normalsize
Thus we have proved that (\ref{ortho13}) is smaller than $C \, \varepsilon$, for any $\varepsilon >0$.

The proof of Theorem \ref{ASYMPT1} is now completed.

\subsection{Proof of Proposition \ref{Feasi1}}

{\color{black}
In the following, for any random function $\beta$, we denote as $c_0(\beta)$ the $0$-th Fourier coefficient of $\beta$, that is, $c_0(\beta)={1\over {2\pi}}\int_0^{2\pi}\beta (t) dt$. Under the conditions on $n,N,M,L$, we prove the convergence in probability
\begin{equation}
\label{V1}
{V}^{(1)}_{n,N,M,L}:=\sum_{|k|\leq L} \bar c_k( \gamma^2_{n,N,M}) \bar c_{-k}( \gamma^2_{n,N,M}) \to c_0(\gamma^4).
\end{equation}

We recall that
\begin{equation}
\label{DueCk}
\bar c_k( \gamma^2_{n,N,M})= c_k( \gamma^2_{n,N,M}) - K 2\pi  c_k(\sigma^4_{n,N,M}),
\end{equation}
according to definitions (\ref{unbiasedK}) and (\ref{biasedK}). Then, we plug (\ref{DueCk}) into (\ref{V1}) and, using the product formula for the Fourier coefficients (see \cite{LMM}), we observe the following convergence in probability: 

$$
\sum_{|k|\leq L} c_k(\gamma^2_{n,N,M})c_{-k}( \gamma^2_{n,N,M}) \,\,\, \to \,\, \, c_0(\gamma^4)+{1\over 9}c_M^4 (1+2\eta({c_N/ \pi}))^2 c_0(\sigma^8)+ {2\over 3}c_M^2(1+2\eta({c_N/ \pi})) c_0(\gamma^2\sigma^4),
$$

$$
 \sum_{|k|\leq L} c_k(\sigma^4_{n,N,M}) c_{-k}(\gamma^2_{n,N,M}) \,\,\, \to \,\, \,   c_0(\gamma^2\sigma^4) +  {1\over 3}c_M^2 (1+2\eta({c_N/ \pi}))c_0(\sigma^8),
$$

$$
 \sum_{|k|\leq L} c_k(\sigma^4_{n,N,M})c_{-k}(\sigma^4_{n,N,M}) \,\,\, \to \,\, \, c_0(\sigma^8).
$$
Then, the convergence in (\ref{V1}) is ensured.
We omit the proof for the convergence of ${  V}^{(2)}_{n,N,M,L}$ and ${  V}^{(3)}_{n,N,M,L}$, which does not contain any novel computation w.r.t. the previous one.}

\subsection{Proof of Theorem \ref{ASYMPT1feasible}}

{\color{black}The theorem follows straightforwardly from the stable convergence in Theorem \ref{ASYMPT1} and Proposition \ref{Feasi1}.}

\subsection{Proof of Theorem \ref{ASYMPT2}}
\label{proofASYMPT2}

The proof relies on the basic decomposition  given in Section \ref{PrelDecom}. First of all, we prove that the bias correction is not needed because
\begin{equation}
\label{correctionNEW}
BB_{M,n,N}^{(i)} +2 BC_{M,n,N}^{(i)}+CC_{M,n,N}^{(i)} = o_p(\rho(n)^{\iota/2}).
\end{equation}

We study the term $BB_{M,n,N}^{(i)}$ defined by (\ref{BBi}).
The term $BC_{M,n,N}^{(i)}$ defined by (\ref{BCi}) and the term $CC_{M,n,N}^{(i)}$ defined by (\ref{CCi}) are analogous to $BB_{M,n,N}^{(i)}$.

Using the It\^o formula, the term $BB_{M,n,N}^{(i)}$ is equal to
\begin{equation}
\label{BIAScomputNEW}
{1\over {M+1}} {1\over {6}} M(M+1)(M+2) {1\over n } \, 2 \, \int_0^{2\pi} \, n\, \int_0^s D_{N,n}^2(s-u) v(u) du \, v(s) ds+o_p(\rho(n)^{1/2}),
\end{equation}
where the order of the martingale part is obtained in Section \ref{biascorrectionterm}.

Now, using Lemma \ref{LemmaDirichlet1}  and noting that $N/n \sim c_N$ and $M/n^{\iota}\sim c_M$,
we have that (\ref{BIAScomputNEW})  has order, in probability, equal to
\begin{equation}
\label{BIAScomput2}
n^{2\iota} {1\over n}.
\end{equation}
It is then enough to observe that $2\iota-1 +\iota/2 <0$,
as long as $\iota < 2/5$. Thus (\ref{correctionNEW}) is proved.

The slower rate of $M$ ensures that the discretization error is still negligible. As for the asymptotic variance, the only term which remains is the bracket
$$\langle \rho(n)^{-\iota/2} 2\, A_M, \rho(n)^{-\iota/2} 2\, A_M \rangle_{2\pi}.$$
Noting that $M/c_M \sim \rho(n)^{-\iota}$, we obtain that its limit in probability is
$$
{1\over {2\pi}} \int_0^{2\pi} {4\over 3} {1\over {c_M}}\gamma^4(t) dt.
$$

\subsection{Proof of Theorem \ref{ASYMPT2feasible}}

{\color{black}
The theorem follows straightforwardly from the stable convergence in Theorem \ref{ASYMPT2} and the convergence in probability of $\Gamma_{n,N,M,L}$ to the asymptotic variance.
The latter is immediately deduced from the following two remarks. First, under the conditions $N \rho(n) \sim c_N$ and $M\rho(n)^{\iota} \sim c_M$, where $\iota \in (0,2/5)$,
then $c_k(\gamma^2_{n,N,M})$, defined in (\ref{biasedK}), converges in probability to $c_k(\gamma^2)$, in virtue of the proof of Theorem \ref{ASYMPT2} and Remark \ref{coefficients}. Secondly, the product formula is applied.}
\section{Appendix B: some properties of the Fej\'{e}r and Dirichlet kernels}
\label{Kernel}

This section resumes some results about the rescaled Dirichlet kernel, defined as
\begin{equation}
\label{dirichlet}
D_N(x):= {1\over {2N+1}} \sum_{|k|\leq N} e^{{\rm i}kx} = {1\over {2N+1}} {\sin ((2N+1)x/2) \over {\sin (x/2)}}
\end{equation}
and the Fej\'{e}r kernel, defined as
\begin{equation}
\label{FejerDef}
F_M(x):= \sum_{|k| \leq M}\left(1-{|k|\over {M+1}}\right) \, e^{{\rm i}kx}={1 \over{M+1}}\left({\sin((M+1) x/2) \over {\sin(x/2)}}\right)^2.
\end{equation}

In the following, we consider a regular partition of the time interval,  maintaining the continuous-time notation used in the proofs of the {CLT}s.

\begin{Lemma}
\label{LemmaFejer}
Suppose that $ {M^{2}/ n} \to a$, as $n,M \to \infty$, for some constant $a > 0$. Then, it holds that:
\begin{equation}
\label{eq:kernel}
\lim_{n, M} \int_{-\pi}^{\pi} F_M(\varphi_n(x))\,dx = \lim_{M} \int_{-\pi}^{\pi} F_M(x)\,dx = 2\pi,
\end{equation}
\begin{equation}
\label{squarefejerTEIC}
\lim_{M,n} \int_{-\pi}^{\pi} {1\over M} F^2_{M}(\varphi_n(x)) dx= \lim_{M} \int_{-\pi}^{\pi} {1\over M} F^2_M(x) dx = {4 \over 3}\pi,
\end{equation}
\begin{equation}
\label{eq:first_derivative_square}
\lim_{n, M} \int_{-\pi}^{\pi}\frac{1}{M^3}|F^{\prime}_M(\varphi_n(x))|^2\,dx = \lim_{M} \int_{-\pi}^{\pi}\frac{1}{M^3}|F^{\prime}_M(x)|^2\,dx = \frac{2}{15}  \pi,
\end{equation}
\begin{equation}
\label{eq:second_derivative_square}
 \lim_{n, M} \int_{-\pi}^{\pi}\frac{1}{M^5} |F^{\prime\prime}_M(\varphi_n(x))|^2\,dx  = \lim_{M} \int_{-\pi}^{\pi} \frac{1}{M^5} |F^{\prime\prime}_M(x)|^2\,dx = \frac{4}{105} \pi.
\end{equation}
Moreover, let
\begin{equation}
K_M(x) := \frac{15 (M + 1) }{M (4 + 6 M + 4 M^2 + M^3)} |F^{\prime}_M(x)|^2,
\end{equation}
\begin{equation}
L_M(x) := \frac{105 (M+1)}{M(2 M^5 + 12 M^4 + 30 M^3 + 40 M^2 + 23 M -2)} |F^{\prime\prime}_M(x)|^2.
\end{equation}
Then, $\left\{K_M(x)\right\}_{M = 1}^{\infty}$ and $\left\{L_M(x)\right\}_{M = 1}^{\infty}$ are families of good kernels.\footnote{See \cite{SteinSha} for the definition.}
\end{Lemma}
\Proof
The results in (\ref{eq:kernel}) and (\ref{squarefejerTEIC}) are proven in \cite{CuTe}, Lemma 5.1. With regards to the first equality in (\ref{eq:first_derivative_square}) and (\ref{eq:second_derivative_square}), it is sufficient to consider the Euler-MacLaurin formula applied to the squared first and second derivative of the Fej\'{e}r kernel. For the sake of completeness, recall that, for a function $f : [-\pi, \pi] \rightarrow \mathbb{R}$ of class $C^{2 p + 1}$, it holds that
\begin{equation}
\label{eq:EuleroMacLaurin}
\int_{[-\pi, \pi]} f(x) \,dx  - \frac{2\pi}{n}\left(\frac{f(\pi) + f(-\pi)}{2} + \sum_{j = 1}^{n-1} f\left(a + j \frac{2\pi}{n}\right)\right)
\end{equation}
$$
= \sum_{k = 1}^{p}\left(\frac{2\pi}{n}\right)^{2 k} \frac{B_{2 k}}{(2 k)!} \left(f^{(2 k - 1)}(-\pi) - f^{(2k-1)}(\pi)\right) + R_{p,n,f},
$$
where $B_{2 k}$ is the $(2 k)$-$th$ Bernoulli number and the rest $R_{p, n, f}$ satisfies
\begin{equation}
|R_{p, n, f}| \leq C_p \left(\frac{2\pi}{n}\right)^{2 p + 1} \int_{-\pi}^{\pi} |f^{(2 p + 1)}(x)|\,dx,
\end{equation}
with $C_p$ a constant depending only on $p$. In particular, let us consider positive integers $k$ and $h$; then, we have that
\begin{equation}
(F_M^{(k)}(x))^{h} = \sum_{|j_1|,\ldots,|j_h| \leq M} \left(1-\frac{|j_1|}{M+1}\right)\cdots\left(1-\frac{|j_h|}{M+1}\right) (\textrm{i} j_1)^k\ldots (\textrm{i} j_h)^{k} e^{\textrm{i} (j_1 +  \ldots + j_h) x}.
\end{equation}
By observing that the number of terms in the summation is $(2M +1)^{h}$, that $|j_1 + \ldots +j_h| \leq h M$, and using the bounds $\left(1-\frac{|j_1|}{M+1}\right)\ldots\left(1-\frac{|j_h|}{M+1}\right) \leq 1$
and $|e^{\textrm{i} (j_1 + \ldots + j_h) x}| \leq 1$,
we have that $|{f}^{(2 p + 1)}| \leq 2^h h^{2 p + 1} M^{(k+1) h + p}$, where
${f} := (F_M^{(k)})^h$.
It follows that $|R_{p, n f} | \leq C_p (2\pi)^{2 p + 2} 2^h h^{2 p + 1} M^{(k+1) h + p} / n^{2 p + 1}$. As a consequence, for $k$ and $h$ fixed and $M^2 / n \rightarrow a$, we have that $|R_{p, n, f}| = O(M^{(k+1) h - p - 1})$. Therefore, for both $f = (1/M^3) |F_M^{(1)}|^{2}$ and $f = (1/M^5) |F_M^{(2)}|^{2}$, we have that $|R_{p, n, f}| = O(M^{- p})$.

We then show  that both the second term on the left hand side  and the first on the right hand side of (\ref{eq:EuleroMacLaurin}) converge  to zero if ${M^2}/{n} \to a$  as $n, M \rightarrow \infty$. For the result of interest, it is sufficient to consider $p = 1$.  First, when $f := (1/M^3) |F_{M}^{(1)}|^2$, the two terms are equal to zero since $F_{M}^{(1)}(\pi) = F_{M}^{(1)}(-\pi) = 0$. Instead, when $f := (1/M^5)|F_{M}^{(2)}|^2$, the term on the left hand side has order $O(M^{-3})$, whereas the term on the right hand side is equal to zero, since $F_M^{(3)}(\pi) = F_M^{(3)}(-\pi) = 0$.

The assertion in (\ref{eq:first_derivative_square}) follows directly from the following calculation
\begin{equation}
\int_{-\pi}^{\pi} \frac{1}{M^3}|F^{\prime}_M(x)|^2 \,dx
= 2\pi \frac{1}{M^3} \sum_{|k| \leq M} \left(1-\frac{|k|}{M+1}\right)^2 k^2 = 2\pi  \frac{M^3 + 4 M^2 + 6 M + 4}{15 M^2 (M + 1)} \rightarrow \frac{2}{15}\pi.
\end{equation}
Similarly, in relation to the assertion in (\ref{eq:second_derivative_square}), one obtains that
$$
\int_{-\pi}^{\pi} \frac{1}{M^5}|F^{\prime\prime}_M(x)|^2 \,dx = 2\pi  \frac{1}{M^5} \sum_{|k| < M} \left(1-\frac{|k|}{M+1}\right)^2 k^4
$$
$$
= 2\pi \frac{2 M^5 + 12 M^4 + 30 M^3 + 40 M^2 + 23 M - 2}{105 M^4 (M+1)} \rightarrow \frac{4}{105} \pi.
$$

It remains to prove that $\left\{K_M(x)\right\}_{M = 1}^{\infty}$ and $\left\{L_M(x)\right\}_{M = 1}^{\infty}$ are families of good kernels. First, consider $K_M$. We observe that $K_M(x) \geq 0$. Then, by using the previous computation, it is easy to show that
$$
\frac{1}{2\pi} \int_{-\pi}^{\pi} K_M(x)\,dx = 1.
$$
Finally, by using the explicit expressions in terms of sine and cosine and the fact that
$|\sin((M+1)x/2)| \leq C (M+1) |x|$ and $|\sin(x/2)| \geq c |x|$ for $|x| \leq \pi$, with $C,c > 0$ suitable constants, we have that
$$
\int_{\delta \leq |x| \leq \pi} |K_M(x)|\,dx = \int_{\delta \leq |x| \leq \pi}   \frac{15 (M+1)}{M (4 + 6 M + 4 M^2 + M^3)}  |F^{\prime}_M(x)|^2\,dx
$$
$$
\leq \frac{15 (M+1)}{M (4 + 6 M + 4 M^2 + M^3)}  \int_{\delta \leq |x| \leq \pi} (C x^{-2})^2 dx
$$
$$
\leq \frac{15 (M+1)}{M (4 + 6 M + 4 M^2 + M^3)} C^2 \int_{\delta \leq |x| \leq \pi}  x^{-4} dx \sim \frac{15 (M+1)}{M (4 + 6 M + 4 M^2 + M^3)} C^2 {1\over {\delta^3}},
$$
which goes to zero as $M \to \infty$.

Analogously, we prove that $\left\{L_M(x)\right\}_{M = 1}^{\infty}$ is a family of good kernels. First, note that $L_M(x) \geq 0$ and
$
\frac{1}{2\pi} \int_{-\pi}^{\pi} L_M(x)\,dx = 1.
$
Moreover, for $A,B$ suitable constants, we have that
$$
\int_{\delta \leq |x| \leq \pi} |L_M(x)|\,dx = \int_{\delta \leq |x| \leq \pi}   \frac{105 (M+1)}{M (-2 + 23 M + 40 M^2 + 30 M^3 + 12 M^4 + 2 M^5)}  (F^{\prime\prime}_M(x))^2\,dx
$$
$$
\leq C^2 \frac{105 (M+1)}{M (-2 + 23 M + 40 M^2 + 30 M^3 + 12 M^4 + 2 M^5)} ( A(M+1) + B (M+1)^{-1})^2  \int_{\delta \leq |x| \leq \pi} x^{-4}\,dx
$$
$$
\sim \frac{105 (M+1)^3 }{M (-2 + 23 M + 40 M^2 + 30 M^3 + 12 M^4 + 2 M^5)}\frac{C^2}{\delta^3},
$$
which converges to zero as $M \rightarrow \infty$.

\begin{Lemma}
\label{LemmaDirichlet1}
\par\noindent
Under the condition $ N/n \to a > 0$, the following results hold.
\par\noindent
i) For any $p>1$, there exists a constant $C_p$ such that
$$\lim_{n,N}\, n \sup_{x\in [0,2\pi]} \int_0^{2\pi} |D_{N}(\varphi_n(x)-\varphi_n(y))|^p dy \leq C_p.$$
ii) It holds that
$$
\lim_{n,N}n \, \int_0^{x} D^2_{N}(\varphi_n(x)-\varphi_n(y)) dy = \pi(1+2\eta(2a))
$$
and, for any $\alpha$-H\"older continuous function $f$, with $\alpha \in (0,1]$,
\begin{equation}
\label{KEY}
\lim_{N,n}n \int_0^{x} D^2_{N}(\varphi_n(x)-\varphi_n(y)) f(y) dy = \pi (1+2\eta (2a)) \, f(x),
\end{equation}
where
\begin{equation}
\label{etaci}
\eta(a):= \frac{1}{2a^2}r(a)(1-r(a)),
\end{equation}
being $r(a)=a-[a]$, with $[a]$ the integer part of $a$.
\par\noindent
iii) For any $\varepsilon >0$,
$$ \lim_{N,n}n \int_0^{x-\varepsilon} |D_{N}(\varphi_n(x)-\varphi_n(y))|^2 dy = 0.$$
\end{Lemma}
Proof. See \cite{ClGl}, Lemma 1 and Lemma 4.

\end{document}